\documentclass[12pt]{article}
\usepackage{amsmath,amssymb,amsthm,latexsym,amscd}

\title{Distributions of countable models \\ of theories with continuum many types\footnote{{\em
Mathematics Subject Classification:} 03C15. \newline\indent \ \ \
The work is supported by RFBR (grant 12-01-00460-a).}}
\author{Roman A. Popkov\footnote{r-popkov@yandex.ru} \ and Sergey V.
Sudoplatov\footnote{sudoplat@math.nsc.ru}}
\date{October 15, 2012}
\begin{document}
\maketitle

\begin{abstract}
We present distributions of countable models and correspondent
structural characteristics of complete theories with continuum
many types: for prime models over finite sets relative to
Rudin--Keisler preorders, for limit models over types and over
sequences of types, and for other countable models of theory.

{\bf Key words:} countable model, theory with continuum many
types, Rudin--Keisler preorder, prime model, limit model, premodel
set.
\end{abstract}

Denote by $\mathcal{T}_c$\index{$\mathcal{T}_c$} the class of all
countable complete, non-small theories $T$, i.~e., of theories
with continual sets $S(T)$ of types. Below, unless otherwise
stated, we shall assume that all theories, under consideration,
belong to the class ${\cal T}_c$ and these theories will be called
{\em unsmall}\index{Theory!unsmall} or {\em theories with
continuum many types}\index{Theory!with continuum many types}.

In general case, for theories in ${\cal T}_c$, there is no
correspondence between types and prime models over tuples that we
observe for small theories (for given theory in ${\cal T}_c$, some
prime models over realizations of types may not exist). Besides
there are continuum many pairwise non-isomorphic countable models
for each of these theories. However as we shall show, in this
case, the structural links for types allow to distribute and to
count the number of prime over finite sets, limit, and other
countable models of a theory like small theories \cite{SuLP, Su08}
and arbitrary countable theories of unary predicates \cite{Pop}.

\bigskip
\centerline{\bf 1. Examples}

\medskip
Recall some basic examples of theories with continuum many types
\cite{KeCh, Spr}:

\medskip
(1) the theory ${\rm Th}(\langle\mathbb N;+,\cdot\rangle)$ of the
standard model of arithmetic on naturals (for any subset $A$ of
the set $P$ of all prime numbers, the set $\Phi(x)$ of formulas
describing the divisibility of an element by a number in $A$ and
its non-divisibility by each number in $P\setminus A$ is
consistent);

\medskip
(2) the theory ${\rm Th}(\langle\mathbb Z;+,0\rangle)$ (there are
continuum many $1$-types by the same reason as in the previous
example);

\medskip
(3) the theory ${\rm Th}(\langle\mathbb Q;+,\cdot,\leq\rangle)$ of
ordered fields (there are $2^\omega$ cuts for the set of
rationals);

\medskip
(4) the theory $T_{\rm sdup}$\index{$T_{\rm sdup}$} of a countable
set of {\em sequentially divisible unary
predicates}\index{Theory!of sequentially divisible!unary
predicates} $S^{(1)}_\delta$, $\delta\in 2^{<\omega}$, with the
following axioms:

$$ \exists^{\geqslant \omega}x \left(S_{\overline{\delta}}(x)\wedge \neg S_{\overline{\delta} \textrm{\^{}}0}(x)
\wedge \neg S_{\overline{\delta} \textrm{\^{}}1}(x\right));$$
$$ S_{\overline{\delta} \textrm{\^{}} {\varepsilon}}(x) \rightarrow
S_{\overline{\delta}}(x),\,\,\varepsilon\in\{0,1\};$$
$$ \neg \exists x (S_{\overline{\delta} \textrm{\^{}} {0}}(x) \wedge  S_{\overline{\delta} \textrm{\^{}} {1}}(x));$$

\medskip
(5) the theory $T_{\rm iup}$\index{$T_{\rm iup}$} of a countable
set of {\em independent unary predicates}\index{Theory!of
independent unary predicates} $P^{(1)}_k$, $k\in\omega$,
axiomatizable by formulas:
$$
\exists x\,(P_{i_1}(x)\wedge\ldots\wedge P_{i_m}(x)\wedge\neg
P_{j_1}(x)\wedge\ldots\wedge\neg P_{j_n}(x)),$$
$\{i_1,\ldots,i_m\}\cap\{j_1,\ldots,j_n\}=\varnothing$ (one get
continuum many $1$-types by consistency of any set of formulas
$\{P^{\delta(k)}_k(x)\mid k\in\omega\}$, $\delta\in 2^\omega$);

\medskip
(6)\footnote{The example is proposed by
E.~A.~Palyutin\index{Palyutin E. A.}.} the theory $T_{\rm
ersiup}$\index{$T_{\rm ersiup}$} of a countable set of {\em
sequentially independent unary predicates $P^{(1)}_k$,
$k\in\omega$, with an equivalence relation}\index{Theory!of
sequentially independent unary predicates with equivalence
relation} $E^{(2)}$, defined by the following axioms:

(a) there are infinitely many $E$-classes and each $E$-class is
infinite;

(b) for any $k\in\omega$ there is unique $E$-class $X_k$
containing infinitely many solutions of each formula
$P^{\delta_0}_{0}(x)\wedge\ldots\wedge P^{\delta_k}_{k}(x)$,
$\delta_0,\ldots,\delta_k\in\{0,1\}$, and $X$ is disjoint with
relations $P_i$, $i>k$; there is a prime model consisting of
$E$-classes $X_k$, $k\in\omega$;

one get continuum many $1$-types in $E$-classes having nonempty
intersections with each predicate $P_k$, $k\in\omega$;

\medskip
(7) the theory $T_{\rm sier}$\index{$T_{\rm sier}$} of a countable
set of {\em sequentially independent equivalence
relations}\index{Theory!of sequentially independent equivalence
relations} $E^{(2)}_n$, $n\in\omega$, with the following axioms:

(a) $\vdash E_{n+1}(x,y)\to E_0(x,y)$, $n\in\omega$;

(b) $\models\forall x,y(E_0(x,y)\to\exists z(E_m(x,z)\wedge
E_n(z,y)))$, $m\ne n$;

(c) each $E_0$-class is infinite and each $E_{n+1}$-class is a
singleton or infinite, $n\in\omega$;

(d) if an $E_{n+1}$-class $X$ is contained in an $E_0$-class $Y$
then $Y$ consists of infinitely many $E_{n+1}$-classes, each of
which is a singleton or infinite, $n\in\omega$;

(e) if $X_{n+1}$ is an infinite $E_{n+1}$-class contained in an
$E_0$-class $Y$ then $Y$ is represented as a union of infinite
intersections $X_1\cap\ldots\cap X_n\cap X_{n+1}$ for
$E_i$-classes $X_i$, $1\leq i\leq n$; moreover, for any
$\delta_i\in\{0,1\}$ the sets $X^{\delta_1}_1\cap\ldots\cap
X^{\delta_n}_n\cap X^{\delta_{n+1}}_{n+1}\cap Y$ are infinite,
$n\in\omega$;

(f) for any $n\in\omega$ there is unique $E_0$-class containing
infinite $E_1\mbox{-},\ldots,E_n$-class and one-element
$E_m$-classes, $n<m$; there is a prime model consisting of these
$E_0$-classes;

there are continuum many $2$-types in $E_0$-classes containing
infinite $E_{n+1}$-classes, $n\in\omega$.

\medskip
The structures $\langle\mathbb N;+,\cdot\rangle$ and
$\langle\mathbb Q;+,\cdot,\leq\rangle$ are prime (since the
universes of there structures equal to ${\rm dcl}(\varnothing)$),
the structure $\langle\mathbb Z;+,0\rangle$ is prime over each its
nonzero element (but it is not prime over $\varnothing$).

The theory $T_{\rm sdup}$ has a prime model and  this model omits
the type $p_\infty(x)$ deduced from the set of formulas describing
the unbounded divisibility of $S_{\overline{\delta}}(x)$ by
$S_{\overline{\delta}\textrm{\^{}} {\varepsilon}}(x)$. Moreover,
the theory $T_{\rm sdup}$ has a prime model over every finite set,
whence there are continuum many pairwise non-isomorphic prime
models over tuples.

The theory $T_{\rm iup}$ does not have prime models over finite
sets. The theories $T_{\rm ersiup}$ and $T_{\rm sier}$, having
prime models over empty set, do not have prime models over
non-principal types.

\bigskip
\centerline{\bf 2. Rudin--Keisler preorders}

\medskip
Consider a theory $T\in\mathcal{T}_c$, a type $p\in S(T)$ and its
realization $\bar{a}$. It is known that all prime models over
realizations of~$p$ are isomorphic. So if there is a {\em prime
model} ${\cal M}(\bar{a})$\index{${\cal M}(\bar{a})$} {\em over
the tuple}~$\bar{a}$,\index{Model!prime over a tuple}, this model
will be usually denoted by ${\cal M}_p$\index{${\cal M}_p$}.

A consistent formula $\varphi(\bar{x})$ of $T$ belonging to an
isolated type in $S(T)$ is called an {\em
$i$-formula},\index{$i$-formula} and if $\varphi(\bar{x})$ does
not belong to isolated types in $S(T)$ then $\varphi(\bar{x})$ is
a {\em ${\rm ni}$-formula}\index{${\rm ni}$-formula}.

Recall \cite{Va} that the prime model of $T$ exists if and only if
every formula being consistent with $T$ is an $i$-formula.

Note that an expansion of any countable structure ${\cal M}$ by
constants for each element transforms this structures to a prime
one. Hence the property of absence of a prime model for a theory
is not preserved under expansions of a theory. Clearly, this
property is not also preserved under restrictions of a theory.

Let $p=p(\bar{x})$ and $q=q(\bar{y})$ be types in $S(T)$.
Following \cite{SuLP, Su041, Su111} we say that $p$ \emph{is
dominated by a type}\index{Type!dominated} $q$, or $p$ \emph{does
not exceed $q$ under the Rudin--Keisler preorder}\index{Type!not
exceed $q$ under the Rudin--Keisler
preorder}\index{Preorder!Rudin--Keisler} (written $p\leq_{\rm RK}
q$\index{$p\leq_{\rm RK}q$}\index{$\leq_{\rm RK}$}), if any model
${\cal M}\models T$ realizing $q$ realizes $p$ too.

By Omitting Types Theorem the condition  $p\leq_{\rm RK}q$ can be
syntactically characterized by the following: there is a {\em
$(q,p)$-formula}\index{$(q,p)$-formula}, i.~e., a formula
$\varphi(\bar{x},\bar{y})$ such that the set
$q(\bar{y})\cup\{\varphi(\bar{x},\bar{y})\}$ is consistent and
$q(\bar{y})\cup\{\varphi(\bar{x},\bar{y})\}\vdash p(\bar{x})$.
Herewith, in contrast small theories, a principal formula
$\varphi(\bar{x},\bar{b})$ with the conditions specified, where
$\models q(\bar{b})$, may not exist. If a principal formula
$\varphi(\bar{x},\bar{b})$ of that form exists, the
$(q,p)$-formula $\varphi(\bar{x},\bar{y})$ is called
\emph{$(q,p)$-principal}\index{Formula!$(q,p)$-principal}.

If $p\leq_{\rm RK}q$ and the models ${\cal M}_p$ and ${\cal M}_q$
exist, we say also that ${\cal M}_p$ \emph{is dominated
by}\index{Model!dominated} ${\cal M}_q$, or ${\cal M}_p$
\emph{does not exceed ${\cal M}_q$ under the Rudin--Keisler
preorder}\index{Model!not exceed $q$ under the Rudin--Keisler
preorder}, and write ${\cal M}_p\leq_{\rm RK}{\cal
M}_q$\index{${\cal M}_p\leq_{\rm RK}{\cal M}_q$}.

If the models ${\cal M}_p$ and ${\cal M}_q$ exist, the condition
${\cal M}_p\leq_{\rm RK}{\cal M}_q$ means that ${\cal M}_q\models
p$, i.~e., some copy ${\cal M}'_p$ of ${\cal M}_p$ is an
elementary submodel of ${\cal M}_q$: ${\cal M}'_p\preceq{\cal
M}_q$.

If the model ${\cal M}_q$ exists then the condition $p\leq_{\rm
RK} q$ implies an existence of $(q,p)$-principal formula, but not
vice versa. Clearly, there is a theory $T$ with types $p$ and $q$
such that $p\leq_{\rm RK} q$, there is a $(q,p)$-principal
formula, and the model ${\cal M}_q$ does not exists (it suffices
to take the theory $T_{iup}$ and $1$-types $p$ and $q$ with
$p=q$).

Obviously, no formula $\varphi(\bar{x},\bar{y})$ can not be both a
$(q,p)$-formula and a $(q,p')$-formula for $p\ne p'$. At the same
time, a fixed formula can be a $(q,p)$-formula even for continuum
many types $q$.

A simplest example of that effect is given by an arbitrary
principal formula $\varphi(\bar{x})$ forming a domination for a
correspondent principal type by all types of given theory. Every
(non)principal type $p(x)\in S(T)$ is dominated by an arbitrary
type $q(\bar{y})\in S(T)$ containing the type $p(y_i)$ and it is
witnessed by the formula $(x\approx y_i)$.

The following example illustrates the mechanism of the domination
for a type by continuum many types in a situation different from
the above.

\medskip
{\bf EXAMPLE 2.1.} Consider a disjunctive union of countable unary
predicates $R_0$ and $R_1$ forming a universe of required
structure. Define a coloring ${\rm Col}\mbox{\rm : }R_0\cup
R_1\to\omega\cup\{\infty\}$ with infinitely many elements for each
color in each predicate $R_0, R_1$. Define a bipartite acyclic
directed graph with a relation $Q$ linking parts $R_0$ and $R_1$
and satisfying the following conditions:

\medskip
${\small\bullet}$  every element $a\in R_1$ of color $m\in\omega$
has infinitely many elements $b\in R_0$ of each color $n\geq m$
such that $(a,b)\in Q$ and there are no elements $c\in R_0$ with
$(a,c)\in Q$ and ${\rm Col}(c)<m$;

\medskip
${\small\bullet}$  every element $a\in R_0$ of color $m\in\omega$
has infinitely many elements $b\in R_1$ of each color $n\leq m$
such that $(b,a)\in Q$ and there are no elements $c\in R_1$ with
$(c,a)\in Q$ and ${\rm Col}(c)>m$.

\medskip
By the construction, for $1$-types $p_i$, isolated by sets
$\{R_i(x)\}\cup\{\neg{\rm Col}_n(x)\mid n\in\omega\}$, $i=0,1$, we
have $p_0\leq_{\rm RK}p_1$ (witnessed by the formula $Q(x,y)$) and
$p_1\not\leq_{\rm RK}p_0$.

That structure is denoted by ${\cal M}_{01}$ and its theory by
$T_{01}$. Expand the structure ${\cal M}_{01}$ by independent
unary predicates $P_k$, $k\in\omega$, on each set defined by the
formula $R_1(x)\wedge{\rm Col}_n(x)$, $n\in\omega$, such that the
type $p_0$ preserves the completeness. Then the type $p_1(x)$ has
continuum many completions $q(x)$, each of which dominates the
type $p_0(x)$ by the formula $Q(x,y)$.

A modification of the example with the theory $T_{\rm sdup}$
instead of $T_{\rm uip}$ leads to the theory for which the formula
$Q(x,y)$ produces the domination of the model ${\cal M}_{p_0}$ to
continuum many models ${\cal M}_q$, where all types $q$ are
completions of the type $p_0$ in $S^1(T_{01})$.~$\Box$

\medskip
Types $p$ and $q$ are said to be
\emph{domination-equivalent}\index{Types!domination-equivalent},
\emph{realization-equivalent}\index{Types!realization-equivalent},
\emph{Rudin--Keisler equivalent}\index{Types!Rudin--Keisler
equivalent}, or \emph{${\rm RK}$-equivalent}\index{Types!${\rm
RK}$-equivalent} (written  $p\sim_{\rm RK} q$\index{$p\sim_{\rm
RK} q$})  if~$p\leq_{\rm RK} q$ and $q\leq_{\rm RK} p$. If
$p\sim_{\rm RK} q$ and the models ${\cal M}_p$ and~${\cal M}_q$
exist then ${\cal M}_p$ and ${\cal M}_q$ are also said to be
\emph{domination-equivalent}\index{Models!domination-equivalent},
\emph{Rudin--Keisler equivalent}\index{Models!Rudin--Keisler
equivalent}, or \emph{${\rm RK}$-equivalent}\index{Models!${\rm
RK}$-equivalent} (written  ${\cal M}_p\sim_{\rm RK}{\cal
M}_q$\index{${\cal M}_p\sim_{\rm RK}{\cal M}_q$}).

As in \cite{Ta4}, types $p$ and $q$ are said to be \emph{strongly
domination-equivalent}\index{Types!domination-equivalent!strongly},
\emph{strongly
realization-equivalent}\index{Types!realization-equivalent!strongly},
\emph{strongly Rudin--Keisler
equivalent}\index{Types!Rudin--Keisler equivalent!strongly}, or
\emph{strongly ${\rm RK}$-equivalent}\index{Types!${\rm
RK}$-equivalent!strongly} (written \ $p\equiv_{\rm
RK}q$\index{$p\equiv_{\rm RK} q$}) if for some realizations
$\bar{a}$ and $\bar{b}$ of $p$ and $q$ accordingly both ${\rm
tp}(\bar{b}/\bar{a})$ and ${\rm tp}(\bar{a}/\bar{b})$ are
principal. Moreover, If the models ${\cal M}_p$ and ${\cal M}_q$
exist, they are said to be \emph{strongly
domination-equivalent}\index{Models!domination-equivalent!strongly},
\emph{strongly Rudin--Keisler
equivalent}\index{Models!Rudin--Keisler equivalent}, or
\emph{strongly ${\rm RK}$-equivalent}\index{Models!${\rm
RK}$-equivalent!strongly} (written ${\cal M}_p\equiv_{\rm RK}{\cal
M}_q$\index{${\cal M}_p\equiv_{\rm RK}{\cal M}_q$}).

Clearly, domination relations form preorders, (strong)
domination-equivalence relations are equivalence relations, and
$p\equiv_{\rm RK}q$ implies $p\sim_{\rm RK}q$.

If ${\cal M}_p$ and ${\cal M}_q$ are not domina\-tion-equivalent
then they are non-isomor\-phic. Moreover, non-isomorphic models
may be found among domination-equivalent ones.

Repeating the proof \cite[Proposition 1]{Su111} we get a syntactic
characterization for an isomorphism of models ${\cal M}_p$ and
${\cal M}_q$. It asserts, as for small theories, that an existence
of isomorphism between ${\cal M}_p$ and ${\cal M}_q$ is equivalent
to the strong domination-equivalence of these models.

\medskip
{\bf PROPOSITION 2.1.} {\em For any types $p(\bar{x})$ and
$q(\bar{y})$ of a theory $T$ having the models ${\cal M}_p$ and
${\cal M}_q$, the following conditions are equivalent:

{\rm (1)} models ${\cal M}_p$ and ${\cal M}_q$ are isomorphic;

{\rm (2)} models ${\cal M}_p$ and ${\cal M}_q$ are strongly
domination-equivalent;

{\rm (3)} there exist $(p,q)$- and $(q,p)$-principal formulas
$\varphi_{p,q}(\bar{y},\bar{x})$ and
$\varphi_{q,p}(\bar{x},\bar{y})$ respectively, such that the set
$$
p(\bar{x})\cup q(\bar{y})\cup\{\varphi_{p,q}(\bar{y},\bar{x}),
\varphi_{q,p}(\bar{x},\bar{y})\}$$ is consistent;

{\rm (4)} there exists a $(p,q)$- and $(q,p)$-principal formula
$\varphi(\bar{x},\bar{y})$, such that the set
$$
p(\bar{x})\cup q(\bar{y})\cup\{\varphi(\bar{x},\bar{y})\}$$ is
consistent.}

\medskip
Denote by ${\rm RK}(T)$\index{${\rm RK}(T)$} the set ${\bf
P}$\index{${\bf P}$} of isomorphism types of models ${\cal M}_p$,
$p\in S(T)$, on which the relation of domination is induced by
$\leq_{\rm RK}$ for models ${\cal M}_p$: ${\rm RK}(T)=\langle{\bf
P};\leq_{\rm RK}\rangle$. We say that isomorphism types ${\bf
M}_1,{\bf M}_2\in{\bf P}$ are
\emph{domination-equivalent}\index{Isomorphism types
domination-equivalent} (written ${\bf M}_1\sim_{\rm RK}{\bf
M}_2$\index{${\bf M}_1\sim_{\rm RK}{\bf M}_2$}) if so are their
representatives.

We consider also the relation $\leq_{\rm RK}$, being defined on
the set $S(T)$ of complete types of a theory $T$. Denote the
structure $\langle S(T);\leq_{\rm RK}\rangle$ by ${\rm
RKT}(T)$\index{${\rm RKT}(T)$}.

Below we investigate links and properties of preordered sets ${\rm
RK}(T)$ and ${\rm RKT}(T)$ as well as links of arbitrary countable
models of a theory with continuum many types.

The following assertion proposes criteria for the existence of the
least element in ${\rm RK}(T)$.

\medskip
{\bf THEOREM 2.2.} {\em For a countable complete theory $T$, the
following conditions are equivalent:

{\rm (1)} the theory $T$ has a prime model;

{\rm (2)} the theory $T$ does not have ${\rm ni}$-formulas;

{\rm (3)} the structure ${\rm RKT}(T)$ has the least $\sim_{\rm
RK}$-class, this class consists of isolated types of $T$ and has a
nonempty intersection with any nonempty set
$[\varphi(\bar{x})]\rightleftharpoons\{p(\bar{x})\in
S(T)\mid\varphi(\bar{x})\in p(\bar{x})\}$.}

\medskip
PROOF. The equivalence $(1)\Leftrightarrow(2)$ forms a criterion
for the existence of prime model of a theory \cite{Va}. The
implications $(1)\Rightarrow(3)$ and $(3)\Rightarrow(2)$ are
obvious.~$\Box$

\medskip
Since theories with continuum many types may not have prime models
over tuples, the limits models may not exist too. Nevertheless the
links between countable models can be observed by the following
generalization of Rudin--Keisler preorder on isomorphism types of
countable models that will be also denoted by $\leq_{\rm
RK}$.\index{$\leq_{\rm RK}$} This generalization extends the
preorder $\leq_{\rm RK}$ for isomorphism types of prime models
over tuples and is based on the inclusion relation for finite
diagrams ${\rm FD}(\mathcal{M})$.

Let ${\bf M}_1$ and ${\bf M}_2$ be isomorphism types of models
$\mathcal{M}_1$ and $\mathcal{M}_2$ (of $T$) respectively. We say
that ${\bf M}_1$ \emph{is dominated by}\index{Type!of
isomorphism!dominated} ${\bf M}_2$ and write ${\bf M}_1\leq_{\rm
RK}{\bf M}_1$\index{${\bf M}_p\leq_{\rm RK}{\bf M}_q$} if each
type in $S^1(\varnothing)$, being realized in ${\bf M}_1$, is
realized in ${\bf M}_2$: ${\rm FD}(\mathcal{M}_1)\subseteq{\rm
FD}(\mathcal{M}_2)$.

Since the relation $\leq_{\rm RK}$ does not depend on
representatives $\mathcal{M}_1$ and $\mathcal{M}_2$ of isomorphism
types ${\bf M}_1$ and ${\bf M}_2$, we shall also write
$\mathcal{M}_1\leq_{\rm
RK}\mathcal{M}_2$\index{$\mathcal{M}_1\leq_{\rm RK}\mathcal{M}_2$}
for the representatives $\mathcal{M}_1$ and $\mathcal{M}_2$ if
${\bf M}_1\leq_{\rm RK}{\bf M}_2$.

We denote by ${\rm CM}(T)$\index{${\rm CM}(T)$} the set ${\bf CM}$
of isomorphism types of countable models of $T$, equipped with the
preorder $\leq_{\rm RK}$ of domination on this set: ${\rm
CM}(T)=\langle{\bf CM};\leq_{\rm RK}\rangle$.

Clearly, ${\rm RK}(T)\subseteq{\rm CM}(T)$. Since having
non-principal types of a countable theory, there is a model of
this theory being not represented in ${\rm RK}(T)$, the equality
${\rm RK}(T)={\rm CM}(T)$ is equivalent to the
$\omega$-categoricity of $T$.

By the definition, a prime model over a type and a limit model
over that type, being non-isomorphic, are domination-equivalent.
Whence any limit models over a common type are also
domination-equivalent.

The generalized relation of domination leads to a classification
of countable models of an arbitrary theory of unary
predicates~\cite{Pop}.

As we pointed out, a series of examples shows that, unlike small
theories, for theories with continuum many types the relations of
domination may not induce least elements (being isomorphism types
of prime models). Besides, by the following example, isomorphism
types of prime models over tuples can quite freely alternate with
the other isomorphism types of countable models.

\medskip
{\bf EXAMPLE 2.2.} We consider a disjunctive union of countable
unary predicates $R_0$ and $R_1$ forming the universe of required
structure. We define a coloring ${\rm Col}\mbox{\rm :
}R_0\to\omega\cup\{\infty\}$ with infinitely many elements for
each color. On the set $R_1$, we put a structure of independent
unary predicates $P_k$, $k\in\omega$. We denote by $T_0$ the
complete theory of the described structure.

Now we fix a dense (in the natural topology) set $X=\{q_m\mid
m\in\omega\}$ of $1$-types containing the formula $R_1(y)$. Using
binary predicates $Q_m$, $m\in\omega$, the type $p_\infty(x)$,
being isolated by the set $\{R_0(x)\wedge\neg{\rm Col}_n(x)\mid
n\in\omega\}$, and neighbourhoods
$R_0(x)\wedge\bigwedge\limits_{i=0}^n\neg{\rm Col}_i(x)$ of
$p_\infty(x)$, we get, in the expanded language, that all types in
$X$ are approximated so that, if the type $p_\infty(x)$ is
realized in a model ${\cal M}$ of expanded theory, then the type
$q_m(y)$ is realized in ${\cal M}$ by the principal formula
$Q_m(a,y)$, where $\models p_\infty(a)$ and $Q_m(a,y)\vdash
q_m(y)$, $m\in\omega$, and the realizability in a model of some
types in $X$ does not imply the realizability of $p_\infty(x)$ in
that model. Thus, a prime model over $p_\infty$ dominates a prime
model over a set $A$, where $A$ consists of realizations of types
in $X$ (one realization of each type).

In turn, the model ${\cal M}_{p_\infty}$ is dominated by a
countable model (being not prime over tuples) which contains a
realization of $p_\infty$ (with realizations of types in $X$) and
a countable set of realizations of $1$-types consistent with
$R_1(x)$ and not belonging to~$X$.~$\Box$

\medskip
Using the notion of dense set of types for the theory $T_{\rm
iup}$ (without the predicate $R_1$) one describes (see \cite{Pop})
the preordered, with respect to $\leq_{\rm RK}$, set ${\bf M}$ of
isomorphism types of countable models of $T_{\rm iup}$. Each
countable model is defined by some countable set of realizations
of a dense set. A model ${\cal M}_1$ is dominated by a model
${\cal M}_2$ if and only if each $1$-type $p$, realized in ${\cal
M}_1$, is realized in ${\cal M}_2$ and the number of realizations
of $p$ in ${\cal M}_1$ does not exceed the number of realizations
of $p$ in ${\cal M}_2$. Since the density of set of types is
preserved under arbitrary removing or adding of a $1$-type, the
set ${\bf M}$ does not have minimal and maximal elements.

\medskip
Example 2.2 illustrates that the absence of prime model of a
theory can be combined with the presence of a prime model over a
tuple. At the same time, as the following proposition asserts,
having a ${\rm ni}$-formula no prime model can not be dominated by
all countable models of theory.

\medskip
{\bf PROPOSITION 2.3.} {\em For any ${\rm ni}$-formula
$\varphi(\bar{x})$ and for any non-principal type $p(\bar{y})\in
S(T)$ there is a non-principal type $q(\bar{x})\in S(T)$
containing the formula $\varphi(\bar{x})$ and which does not
dominate the type $p(\bar{y})$.}

\medskip
PROOF. By Omitting Type Theorem, there is a countable model ${\cal
M}$ of $T$ omitting the type $p(\bar{y})$. At the same time, by
consistency of $\varphi(\bar{x})$ there is a tuple~$\bar{a}$ such
that ${\cal M}\models\varphi(\bar{a})$. The type
$q(\bar{x})\rightleftharpoons{\rm tp}(\bar{a})$, contains the
formula $\varphi(\bar{x})$ and, by the definition, does not
dominate the type $p(\bar{y})$.~$\Box$

\medskip
Since each consistent conjunction of ${\rm ni}$-formula
$\varphi(\bar{x})$ and a formula $\psi(\bar{x})$ is again a ${\rm
ni}$-formula, there are infinitely many types $q(\bar{x})\in S(T)$
containing the formula $\varphi(\bar{x})$ and do not dominating
the type $p(\bar{y})$. Moreover, in a series of examples of $T$
like above, there are uncountably many these types since otherwise
there is a countable expansion $T'$ of $T$ with new predicates
$Q_n(\bar{x},\bar{y})$, $n\in\omega$, producing the isolation of
each type $r(\bar{x})\in S(T')$, containing $\varphi(\bar{x})$, by
its restriction to the language of $T$, and the domination of
$p(\bar{y})$ by each type $q(\bar{x})$. Since the formula
$\varphi(\bar{x})$ is again a ${\rm ni}$-formula, we get a
contradiction by Proposition 2.4.

Note that if a type $p(\bar{y})$ is not dominated by a type
$q(\bar{x})$ then, introducing new independent predicates
$P_k(\bar{x})$, $k\in\omega$, transforming a neighbourhood of
$q(\bar{x})$ to a ${\rm ni}$-formula and $q(\bar{x})$ to
$2^\omega$ completions, we get a theory such that $p(\bar{x})$ is
not dominated by continuum many types. By a similar way, as in
Example 2.1, if a type $p(\bar{y})$ is dominated by a type
$q(\bar{x})$ then, in an expansion, the type $p(\bar{y})$ is
dominated by continuum many completions of $q(\bar{x})$.

Note also that a structure ${\rm RKT}(T)$ can have a minimal but
not least $\sim_{\rm RK}$-class. Indeed, expanding the theory
$T_{\rm iup}$ by binary predicates, one can obtain a dense set $S$
of $1$-types, each of which is domination-equivalent with the
other, and the absence of prime model is preserved (it can be done
by a countable set of new binary predicates, each of which is
responsible for the domination-equivalence of two $1$-types in the
given dense set, and this domination-equivalence is obtained by
approximations for neighbourhoods of given types). The set $S$ and
types, domination-equivalent to types in $S$, form a minimal
$\sim_{\rm RK}$-class. By similar expansions, one get countably
many minimal classes.

Together with Example 2.2 and Proposition 2.4, Example 2.1
illustrate a mechanism of domination of a non-principal type by
all non-principal types of a theory with continuum many types and
without ${\rm ni}$-formulas.

Having the features, in the following section, we propose a list
of some basic properties of structures ${\rm RKT}(T)$ for theories
$T$ in $\mathcal{T}_c$.\footnote{Recall that for countable
structures ${\rm RKT}(T)$ the basic properties (the countable
cardinality, the upward direction, the countability of $\sim_{\rm
RK}$-classes, the presence of the least $\sim_{\rm RK}$-classes)
are presented in \cite{Su111}.}

\bigskip
\centerline{\bf 3. Premodel sets}

\medskip
A {\em height}\index{Height of preordered set} ({\em
width})\index{Width of preordered set} of preordered set $\langle
X;\leq\rangle$ is a supremum of cardinalities for its
$\leq$-(anti)chains consisting of pairwise non-$\sim$-equivalent
elements, where $\sim\,\,\rightleftharpoons(\leq\cap\geq)$. Recall
\cite{SuLP}, that if $a\in X$ then the set
$\bigtriangleup(a)$\index{$\bigtriangleup(a)$} (respectively
$\bigtriangledown(a)$\index{$\bigtriangledown(a)$}) of all
elements $x$ in $X$, for which $x\leq a$ ($a\leq x$), is a {\em
lower {\rm (}upper{\rm )}
cone}\index{Cone!lower}\index{Cone!upper} of $a$.

A continual preordered upward directed set $\langle
X;\leq$~$\rangle$ is called {\em premodel}\index{Set!premodel} if
it has:

\medskip
${\small\bullet}$ countably many elements under each element $a\in
X$ (i.~e., $|\bigtriangleup(a)|=\omega$);

\medskip
${\small\bullet}$ only countable $\sim$-classes (i.~e.,
$|\bigtriangleup(a)\cap\bigtriangledown(a)|=\omega$ for any $a\in
X$);

\medskip
${\small\bullet}$ countable, or continual and coinciding with $X$,
co-countable, or co-continual set of common elements over any
elements $a_1,\ldots,a_n\in X$ (i.~e.,
$|\bigtriangledown(a_1)\cap\ldots\cap\bigtriangledown(a_n)|=\omega$,
or
$|\bigtriangledown(a_1)\cap\ldots\cap\bigtriangledown(a_n)|=2^\omega$
and $\bigtriangledown(a_1)\cap\ldots\cap\bigtriangledown(a_n)=X$,
$|X\setminus(\bigtriangledown(a_1)\cap\ldots\cap\bigtriangledown(a_n))|=\omega$,
or
$|X\setminus(\bigtriangledown(a_1)\cap\ldots\cap\bigtriangledown(a_n))|=2^\omega$);

\medskip
${\small\bullet}$ the countable height.

\medskip
{\bf PROPOSITION 3.1.} {\em If  $|S(T)|=2^\omega$ then the
structure ${\rm RKT}(T)$ is premodel.}

\medskip
PROOF. The structure ${\rm RKT}(T)$ is upward directed since types
$p(\bar{x})$, $q(\bar{y})\in S(T)$, where $\bar{x}$ and $\bar{y}$
are disjoint, are dominated by any type $r(\bar{x},\bar{y})\supset
p(\bar{x})\cup q(\bar{y})$ in $S(T)$.

As $T$ is countable, the set of formulas of $T$ is also countable
and each type dominates at most countably many types. Having
countably many types, being domination-equivalent with a given
type (for instance, a type ${\rm tp}(\bar{a})$ is
domination-equivalent with types ${\rm
tp}(\bar{a}\,\hat{\,}\,\bar{a})$, ${\rm
tp}(\bar{a}\,\hat{\,}\,\bar{a}\,\hat{\,}\,\bar{a}),\ldots$), we
get that any type is domination-equivalent with countably many
types of $T$.

Since each formula witnesses on domination of a type to at most
countably many, or continuum and co-continuum many types, and
there are countably many formulas of $T$, then any types
$p_1,\ldots,p_n$ lay under countably many, or continuum many and
coinciding with $S(T)$, co-countably many, or co-continuum many
types.

As each type dominates countably many types, the height of ${\rm
RKT}(T)$ is at most countable. At the same time the height can not
be finite since its finiteness, the upward direction of ${\rm
RKT}(T)$, and the countable domination imply that ${\rm RKT}(T)$
is countable in spite of $|S(T)|=2^\omega$.~$\Box$

\medskip
Since each $\sim_{\rm RK}$-class of a countable theory $T$ is
countable and each type dominates countably many types, the
ordered factor set ${\rm RKT}(T)/\!\!\sim_{\rm RK}$ can be
linearly ordered only for small $T$. Moreover, as the height of
${\rm RKT}(T)$ is countable for $T\in\mathcal{T}_c$, this
factor-set has continuum many incomparable elements, i.~e., the
width is continual:

\medskip
{\bf PROPOSITION 3.2.} {\em The width of any premodel set $\langle
X;\leq$~$\rangle$ is continual.}

\medskip
PROOF. Assume the contrary that the width of a preordered set
$\langle X;\leq$~$\rangle$ is not continual. Consider a maximal
antichain $Y$. By the assumption, we have $|Y|=\lambda<2^\omega$.
We link each element $y\in Y$ with a maximal chain $C_y$. Each
chain $C_y$ is countable since the height is countable and each
$\sim$-class is countable too. Now we note that
$X=\bigcup\{\bigtriangleup(c)\mid c\in C_y,y\in Y\}$ since
$\langle X;\leq$~$\rangle$ is upward directed. Then, as each lower
cone $\bigtriangleup(c)$ is countable, we obtain
$|X|\leq\lambda\cdot\omega\cdot\omega<2^\omega$ that contradicts
the condition $|X|=2^\omega$.~$\Box$

\bigskip
\centerline{\bf 4. Distributions for countable models of a theory}
\centerline{\bf by $\leq_{\rm RK}$-sequences}

\medskip
Recall that, by Tarski--Vaught criterion\index{Tarski--Vaught
criterion}, a set $A$ in a structure ${\cal M}$ of language
$\Sigma$ forms an elementary substructure if and only if for any
formula $\varphi(x_0,x_1,\ldots,x_n)$ of the language $\Sigma$ and
for any elements $a_1,\ldots,a_n\in A$ if ${\cal M}\models\exists
x_0\,\varphi(x_0,a_1,\ldots,a_n)$ then there is an element $a_0\in
A$ such that ${\cal M}\models\varphi(a_0,a_1,\ldots,a_n)$. It
means that each formula $\varphi(\bar{x})$ over a finite set
$A_0\subseteq A$ and belonging to a type over $A_0$ has a
realization $\bar{a}\in A$.

Let ${\cal M}$ be a model of a countable theory $T$ and ${\bf
q}\rightleftharpoons(q_n)_{n\in\omega}$ be a {\em $\leq_{\rm
RK}$-sequence}\index{$\leq_{\rm RK}$-sequence} of types of $T$,
i.~e., a sequence of non-principal types $q_n$ with $q_n\leq_{\rm
RK}q_{n+1}$, $n\in\omega$. We denote by $U({\cal M},{\bf
q})$\index{$U({\cal M},{\bf q})$} the set of all realizations in
${\cal M}$ of types of $T$, being dominated by some types in ${\bf
q}$. The $\leq_{\rm RK}$-sequence ${\bf q}$ is called {\em
elementary submodel}\index{$\leq_{\rm RK}$-sequence!elementary
submodel} if for any consistent formula $\varphi(\bar{y})$ of $T$
some type in ${\bf q}$ dominates a type $p(\bar{y})\in S(T)$
containing the formula $\varphi(\bar{y})$, and if the formula
$\varphi(\bar{y})$ is equal to $\exists x\,\psi(x,\bar{y})$ then
the type $p(\bar{y})$ is extensible to a type $p'(x,\bar{y})\in
S(T)$ dominated by a type in ${\bf q}$ and such that
$\psi(x,\bar{y})\in p'$.

\medskip
{\bf THEOREM 4.1.} {\em For any $\omega$-homogeneous model ${\cal
M}$ of a countable theory $T$ and for any $\leq_{\rm RK}$-sequence
${\bf q}$ of types in $S(T)$, realized in $\mathcal{M}$, the
following conditions are equivalent:

$(1)$ some {\rm (}countable{\rm )} subset of $U({\cal M},{\bf q})$
is a universe of elementary submodel of ${\cal M}$;

$(2)$ ${\bf q}$ is an elementary submodel $\leq_{\rm
RK}$-sequence.}

\medskip
PROOF. $(1)\Rightarrow (2)$ is implied by Tarski--Vaught
criterion.

$(2)\Rightarrow (1)$. Let ${\bf q}$ be an elementary submodel
$\leq_{\rm RK}$-sequence. Using elements of $U({\cal M},{\bf q})$
we construct, by induction, a countable elementary submodel of
${\cal M}$. On initial step we enumerate, by natural numbers, all
consistent, with $T$, formulas $\varphi(x,\bar{y})$ such that the
enumeration $\nu$ starts with some formula $\varphi_0(x)$ and each
formula has infinitely many numbers. We choose a realization
$a_0\in U({\cal M},{\bf q})$ of the formula $\varphi_0(x)$ and put
$A_0\rightleftharpoons\{a_0\}$. Assume that, on step $n$, a finite
set $A_n\subset U({\cal M},{\bf q})$ is defined, the type of this
set is dominated by some type in ${\bf q}$, and all possible
tuples of elements in $A_n$ are substituted in initially
enumerated formulas $\varphi(x,\bar{y})$ instead of tuples
$\bar{y}$ such that there are infinitely many numbers for each
formula, where tuples of elements in $A_n$ are not substituted. We
assume also that the results
$(\varphi(x,\bar{y}))^{\bar{y}}_{\bar{a}}$ of substitutions have
the same numbers as before, a substitution is carried out for the
formula with the number $n+1$, and this formula has the form
$\varphi(x,\bar{a})$. If ${\cal M}\models\neg\exists
x\varphi(x,\bar{a})$, we put $A_{n+1}\rightleftharpoons A_n$. If
${\cal M}\models\exists x\varphi(x,\bar{a})$, we add fictitiously
to the tuple $\bar{a}$ all missing elements of~$A_n$ and choose an
existing, by conjecture, type $p'(x,\bar{y})$ extending the type
$p(\bar{y})={\rm tp}(A_n)$, where $\varphi(x,\bar{y})\in p'$ and
the types $p$, $p'$ are dominated by some types in ${\bf q}$. We
take for $a_{n+1}$ a realization in~$U({\cal M},{\bf q})$ of the
type $p'(x,A_n)$ (that exists since the model ${\cal M}$ is
$\omega$-homogeneous) and put $A_{n+1}\rightleftharpoons
A_n\cup\{a_{n+1}\}$.

It is easy to see, using a mechanism of consistency \cite{ErPa},
that $\bigcup\limits_{n\in\omega}A_n$ is a universe of required
elementary submodel of ${\cal M}$.~$\Box$

\medskip
Since every $\omega$-saturated structure is $\omega$-homogeneous,
Theorem 4.1 implies

\medskip
{\bf COROLLARY 4.2.} {\em For any $\omega$-saturated model ${\cal
M}$ of a countable theory $T$ and for any $\leq_{\rm RK}$-sequence
${\bf q}$ of types in $S(T)$, the following conditions are
equivalent:

$(1)$ some {\rm (}countable{\rm )} subset of $U({\cal M},{\bf q})$
is a universe of elementary submodel of ${\cal M}$;

$(2)$ ${\bf q}$ is an elementary submodel $\leq_{\rm
RK}$-sequence.}

\medskip
Note that, in the proof of Theorem 4.1, we essentially use that
the model $\mathcal{M}$ is $\omega$-homogeneous and all types of
the sequence ${\bf q}$ are realized in $\mathcal{M}$. Possibly the
types of a $\leq_{\rm RK}$-sequence ${\bf q}$ are not realized in
an $\omega$-homogeneous model $\mathcal{M}$ but are realized in in
some other $\omega$-homogeneous model $\mathcal{M}'$, where
Theorem 4.1 can be applied.

\medskip
{\bf EXAMPLE 4.1.} Consider the theory $T_{\rm iup}$. By Theorem
4.1, each countable model of $T_{\rm iup}$ realizes a dense set
$X$ of $1$-types (where $\bigcup X$ contains all formulas
$$
P_{i_1}(x)\wedge\ldots\wedge P_{i_m}(x)\wedge\neg
P_{j_1}(x)\wedge\ldots\wedge\neg P_{j_n}(x),$$
$\{i_1,\ldots,i_m\}\cap\{j_1,\ldots,j_n\}=\varnothing$) and vice
versa, for each countable dense set $X$ of types, there is an
($\omega$-homogeneous) model of $T_{\rm iup}$ such that the set of
types of elements equals to $X$.

Take two countable disjoint dense sets $P_0$ and $P_1$ of
$1$-types, and $\omega$-homogeneous models $\mathcal{M}_0$ and
$\mathcal{M}_1$ containing exactly one realization of each type in
$P_0$ and $P_1$ respectively. Then there are $\leq_{\rm
RK}$-sequences ${\bf q}_i$ of types with realizations from given
sets of realizations of types in $P_i$, $i=0,1$. Here, all types
in ${\bf q}_i$ are realized $\mathcal{M}_i$ and are omitted in
$\mathcal{M}_{1-i}$, $i=0,1$.~$\Box$

\medskip
By Theorem 4.1, each elementary submodel $\leq_{\rm RK}$-sequence
${\bf q}$ corresponds to some set of isomorphism types of
countable models of a theory $T$, which can vary from $1$ to
$2^\omega$. We denote this set by $I^{m}_{\bf
q}(T)$\index{$I^{m}_{\bf q}(T)$}.

The sets $I^m_{\bf q}(T)$ can have nonempty intersections (for
instance, having a prime model ${\cal M}_0$ its isomorphism type
belongs to each set $I^m_{\bf q}(T)$) and can be disjoint (as in
Example 4.1).

Distributing isomorphism types of countable model to pairwise
disjoint sets, related to $\leq_{\rm RK}$-sequences ${\bf q}$ (and
not related to the other $\leq_{\rm RK}$-sequences) and denoting
the cardinalities of these sets by~$I_{\bf q}$\index{$I_{\bf q}$},
we have the equality
$$I(T,\omega)=\sum\limits_{\bf q}I_{\bf
q}=2^\omega.$$

\medskip
\centerline{\bf 5. Three classes of countable models}
\medskip

Recall \cite{SuLP, Su08} that a model ${\cal M}$ of a theory $T$
is called {\em limit}\index{Model!limit} if ${\cal M}$ is not
prime over tuples and ${\cal M}=\bigcup\limits_{n\in\omega}{\cal
M}_n$ for some elementary chain $({\cal M}_n)_{n\in\omega}$ of
prime models of $T$ over tuples. In this case the model ${\cal M}$
is said to be {\em limit over a sequence ${\bf q}$ of
types}\index{Model!limit!over a sequence of types}, where ${\bf
q}=(q_n)_{n\in\omega}$, ${\cal M}_n={\cal M}_{q_n}$, $n\in\omega$.
If a cofinite subset of the set of types $q_n$ is a singleton
containing a type $p$ then the limit model over ${\bf q}$ is said
to be {\em limit over the type}\index{Model!limit!over type} $p$.

Consider a countable complete theory $T$. Denote by ${\bf
P}\rightleftharpoons{\bf P}(T)$,\index{${\bf P}(T)$}\index{${\bf
P}$} ${\bf L}\rightleftharpoons{\bf L}(T)$\index{${\bf
L}(T)$}\index{${\bf L}$}, and ${\bf NPL}\rightleftharpoons{\bf
NPL}(T)$\index{${\bf NPL}(T)$}\index{${\bf NPL}$} respectively the
set of prime over tuples, limit, and other countable models of
$T$, and by $P(T)$,\index{$P(T)$} $L(T)$\index{$L(T)$}, and ${\rm
NPL}(T)$\index{${\rm NPL}(T)$} the cardinalities of these sets.

By the definition, each value $P(T)$, $L(T)$, and ${\rm NPL}(T)$
may vary from $0$ to $2^\omega$ and the following equality holds:
$$
I(T,\omega)=P(T)+L(T)+{\rm NPL}(T).
$$
Since $I(T,\omega)=2^\omega$ for theories $T$ in $\mathcal{T}_c$,
some value $P(T)$, $L(T)$, or ${\rm NPL}(T)$ is equal to
$2^\omega$.

The tuple $(P(T),L(T),{\rm NPL}(T))$ is called a {\em triple of
distribution of countable models of $T$}\index{Triple of
distribution of countable models of theory} and is denoted by
${\rm cm}_3(T)$.\index{${\rm cm}_3(T)$}

\medskip
{\bf Definition.} A theory $T$ is called {\em
$p$-zero}\index{Theory!$p$-zero} (respectively {\em
$l$-zero}\index{Theory!$l$-zero}, {\em ${\rm
npl}$-zero}\index{Theory!${\rm npl}$-zero}) if $P(T)=0$
(respectively $L(T)=0$, ${\rm NPL}(T)=0$).

A theory $T$ is called {\em
$p$-categorical}\index{Theory!$p$-categorical} (respectively {\em
$l$-categorical}\index{Theory!$l$-categorical}, {\em ${\rm
npl}$-categorical}\index{Theory!${\rm npl}$-categorical}) if
$P(T)=1$ (respectively $L(T)=1$, ${\rm NPL}(T)=1$).

A theory $T$ is called {\em
$p$-Ehrenfeucht}\index{Theory!$p$-Ehrenfeucht} (respectively {\em
$l$-Ehrenfeucht}\index{Theory!$l$-Ehrenfeucht}, {\em ${\rm
npl}$-Ehrenfeucht}\index{Theory!${\rm npl}$-Ehrenfeucht}) if
$1<P(T)<\omega$ (respectively $1<L(T)<\omega$, $1<{\rm
NPL}(T)<\omega$).

A theory $T$ is called {\em
$p$-countable}\index{Theory!$p$-countable} (respectively {\em
$l$-countable}\index{Theory!$l$-countable}, {\em ${\rm
npl}$-countable}\index{Theory!${\rm npl}$-countable}) if
$P(T)=\omega$ (respectively $L(T)=\omega$, ${\rm NPL}(T)=\omega$).

A theory $T$ is called {\em
$p$-continual}\index{Theory!$p$-continual} (respectively {\em
$l$-continual}\index{Theory!$l$-continual}, {\em ${\rm
npl}$-continual}\index{Theory!${\rm npl}$-continual}) if
$P(T)=2^\omega$ (respectively $L(T)=2^\omega$, ${\rm
NPL}(T)=2^\omega$).

\medskip
By the definition, each $p$-zero theory is $l$-zero.

Recall \cite{SuLP, Su08} that the $p$-categoricity of a small
theory $T$ is equivalent to its countable categoricity as well as
to the absence of limit models. The $p$-Ehrenfeuchtness of $T$
means that the structure ${\rm RK}(T)$ is finite and has at least
two elements. The theory $T$ is Ehrenfeucht if and only if $T$ is
$p$-Ehrenfeucht and $L(T)<\omega$. Besides every small theory is
${\rm npl}$-zero, i.~e., each its countable model is prime over a
tuple or is limit. Since by Vaught's and Morley's theorems
\cite{Va, Mo2},
$I(T,\omega)\in(\omega\setminus\{0,2\})\cup\{\omega,\omega_1,2^\omega\}$
and for small theories $T$ the inequalities $1\leq P(T)\leq\omega$
hold, we have the following

\medskip
{\bf THEOREM 5.1.} {\em For any small theory $T$ the triple ${\rm
cm}_3(T)$ has one of the following values:

$1)$ $(1,0,0)$ {\rm (}any $p$-categorical theory, being
$\omega$-categorical, is $l$-zero and ${\rm npl}$-zero{\rm )};

$2)$ $(\lambda_1,\lambda_2,0)$, where $2\leq\lambda_1\leq\omega$,
$\lambda_2\in(\omega\setminus\{0\})\cup\{\omega,\omega_1,2^\omega\}$
{\rm (}for non-$\omega$-categorical small theories{\rm )}.}

\medskip
As shown in \cite{SuLP, Su08}, all values, pointed out in Theorem
5.1 (for $\lambda_2\ne\omega_1$) have realizations in the class of
small theories.

Similarly Theorem 5.1, for the classification of theories in the
class $\mathcal{T}_c$, the problem arises for the description of
all possible triples $(\lambda_1,\lambda_2,\lambda_3)$ realized by
${\rm cm}_3(T)$ for theories $T\in\mathcal{T}_c$.

Examples in Section 1 confirm the realizability of triples
$(0,0,2^\omega)$ and $(2^\omega,2^\omega,0)$ in the class
$\mathcal{T}_c$ (by the $p$-zero, ${\rm npl}$-continual theory
$T_{\rm iup}$ and the $p$-continual, ${\rm npl}$-zero theory
$T_{\rm sdup}$ respectively). Some fusion of theories $T_{\rm
iup}$ and $T_{\rm sdup}$ substantiates the realizability of triple
$(2^\omega,2^\omega,2^\omega)$. E.~A.~Palyutin\index{Palyutin E.
A.} noted that the theory $T_{\rm ersiup}$ realizes the triple
$(1,0,2^\omega)$. This triple is also realized by the theory
$T_{\rm sier}$.

The following theorem produces a characterization for the class of
${\rm npl}$-zero theories.

\medskip
{\bf THEOREM 5.2.} {\em A countable model ${\cal M}$ of a theory
$T\in{\cal T}_c$ is prime over a finite set or limit if and only
if each tuple $\bar{a}\in M$ is extensible to a tuple $\bar{b}\in
M$ such that each consistent formula $\varphi(\bar{x},\bar{b})$ is
an $i$-formula.}

\medskip
PROOF. If for a tuple $\bar{b}\in M$ every consistent formula
$\varphi(\bar{x},\bar{b})$ is an $i$-formula then there is a model
${\cal M}(\bar{b})\preccurlyeq {\cal M}$. Repeating the proof of
Proposition 1.1.7 in \cite{SuLP} or of Proposition 4.1 in
\cite{Su08} in respect that any tuple $\bar{a}$ is extensible to a
tuple $\bar{b}$ of described form, we get a representation of
${\cal M}$ as a union of elementary chain of prime models over
finite sets. Thus, ${\cal M}$ is prime over a finite set or limit.

If a tuple $\bar{a}\in M$ is not extensible to a tuple $\bar{b}\in
M$ such that each consistent formula $\varphi(\bar{x},\bar{b})$ is
$i$-formula, then $\bar{a}$ is not contained in prime models over
tuples, being elementary submodels of ${\cal M}$, whence the model
${\cal M}$ is neither prime over a tuple nor limit.~$\Box$

\medskip
Theorem 5.2 implies

\medskip
{\bf COROLLARY 5.3.} {\em A theory $T\in{\cal T}_c$ is ${\rm
npl}$-zero if and only if for any {\rm (}countable{\rm )} model
${\cal M}$ of $T$ each tuple $\bar{a}\in M$ is extensible to a
tuple $\bar{b}\in M$ such that every consistent formula
$\varphi(\bar{x},\bar{b})$ is an $i$-formula.}

\medskip
Below we describe some families of triples
$(\lambda_1,\lambda_2,\lambda_3)$ that cannot be realized by ${\rm
cm}_3(T)$, where $T\in\mathcal{T}_c$.

\medskip
{\bf PROPOSITION 5.4.} {\em There is no theory $T\in\mathcal{T}_c$
such that ${\rm cm}_3(T)$ has any of the following:

$(1)$ $(\lambda_1,2^\omega,\lambda_3)$, where
$\lambda_1,\lambda_3<2^\omega$;

$(2)$ $(2^\omega,\lambda_2,\lambda_3)$, where
$\lambda_2,\lambda_3<2^\omega$.}

\medskip
PROOF. (1) If $P(T)<2^\omega$ and ${\rm NPL}(T)<2^\omega$ then
there are less than continuum many types that realized in models
representing isomorphism types in the classes ${\bf P}(T)$ and
${\bf NPL}(T)$. Since each type, realized in a limit model, is
also realized in a prime model over a tuple, there are continuum
many types, being not realized in countable models of $T$, that is
impossible.

(2) Assume that ${\rm NPL}(T)<2^\omega$. Then there are
$<2^\omega$ types in $S(T)$, over which prime models do not exist.
Therefore, for any type $p\in S(T)$ there are continuum many types
$q\in S(T)$ extending $p$ and having models $\mathcal{M}_q$. Since
there are continuum many types $q$ and the model $\mathcal{M}_p$
is countable, then there are continuum many these
non-domination-equivalent types $q$ dominating $p$ and not
dominated by $p$. Whence, for any model $\mathcal{M}_p$ there are
continuum many possibilities for elementary extensions by pairwise
non-isomorphic models $\mathcal{M}_q$, being non-isomorphic to
$\mathcal{M}_p$. Since the process of extension of models
$\mathcal{M}_p$ by continuum many models $\mathcal{M}_q$ can be
continued unboundedly many times, there are continuum many
pairwise non-isomorphic limit models, i.~e.,
$L(T)=2^\omega$.~$\Box$

\medskip
The following proposition gives a sufficient condition for the
existence of continuum many prime models over finite sets in the
assumption of uncountably many these models.

\medskip
{\bf PROPOSITION 5.5.} {\em Let there be uncountably many types
$p(\bar{x})$ of a theory $T\in\mathcal{T}_c$ such that for each
formula $\varphi(\bar{a},\bar{y})$, $\models p(\bar{a})$, there is
a principal formula $\psi(\bar{a},\bar{y})$ with
$\psi(\bar{a},\bar{y})\vdash\varphi(\bar{a},\bar{y})$ and this
formula can be chosen independently of the types $p$. Then
$P(T)=2^\omega$.}

\medskip
PROOF. Since there are uncountably many types $p(\bar{x})$, we
have neighbourhoods $\chi_\delta(\bar{x})$ of these types,
$\delta\in 2^{<\omega}$, each of which belongs to uncountably many
given types $p(\bar{x})$ and satisfies the following conditions:

\medskip
${\small\bullet}$
$\chi_\delta(\bar{x})\equiv(\chi_{\delta\,\hat{\,}\,0}(\bar{x})\vee\chi_{\delta\,\hat{\,}\,1}(\bar{x}))$;

\medskip
${\small\bullet}$
$\models\neg\exists\bar{x}(\chi_{\delta\,\hat{\,}\,0}(\bar{x})\wedge\chi_{\delta\,\hat{\,}\,1}(\bar{x}))$.

\medskip
For each sequence $\delta\in 2^\omega$, the local consistency
implies the consistency of the set $\Phi_\delta(\bar{x})$ of
formulas $\chi_{\delta\upharpoonright n}(\bar{x})$, $n\in\omega$.
Whence there are continuum many types in
$S^{l(\bar{x})}(\varnothing)$. Moreover, since the formulas $\psi$
can be chosen by $\varphi$ independently of realizations of types
$p$, by compactness each set $\Phi_\delta(\bar{x})$ has a
completion $q(\bar{x})\in S(\varnothing)$ such that for any
consistent formula $\varphi(\bar{a},\bar{y})$, $\models
q(\bar{a})$, there is a principal formula $\psi(\bar{a},\bar{y})$
with $\psi(\bar{a},\bar{y})\vdash\varphi(\bar{a},\bar{y})$ and
this formula does not depend of $q$ as well as it was independent
of $p$. Thus, there is a model $\mathcal{M}_q$ and there are
continuum many these models, i.~e., $P(T)=2^\omega$.~$\Box$

\medskip
By Proposition 5.5, we have a partial solution of a variant of the
Vaught problem, being formulated by E.~A.~Palyutin\index{Palyutin
E. A.} as the implication $P(T)>\omega\Rightarrow P(T)=2^\omega$.
Namely, this implication is true for prime models over
realizations of types $p$ having the specified, in the
proposition, {\em uniform choice property}\index{Property!of
uniform choice} of formulas $\psi$ by formulas $\varphi$.

\medskip
\centerline{\bf 6. Operators acting on a class of structures}
\medskip

Consider a non-principal $1$-type $p_{\infty}(x)$ and formulas
$\varphi_n(x) ~\in~ p_{\infty}(x)$, $n \in \omega$, such that
$\varphi_0(x) = (x \approx x)$, $\vdash \varphi_{n+1}(x)
\rightarrow \varphi_{n}(x)$, $\{ \varphi_n(x) \mid n \in \omega \}
\vdash p_{\infty}(x)$. The formula ${\rm
Col}_n(x)\rightleftharpoons\varphi_n(x) \wedge \neg
\varphi_{n+1}(x)$\index{${\rm Col}_n(x)$} is the $n$-\emph{th
approximation}\index{Approximation!$n$-th} of $p_{\infty}(x)$ or
the $n$-\emph{th  color}.\index{Color!$n$-th} Then the type
$p_{\infty}(x)$ is isolated by the set $\{\neg {\rm Col}_n \mid n
\in \omega \}$ of formulas.

\medskip

The \emph{operator of continual partition}\index{Operator!of
continual partition} ${\rm icp}(\mathcal{A}, \mathcal{A}_0, Y,
\{R_i^{(2)}\}_{i \in \omega})$\index{${\rm icp}(\mathcal{A},
\mathcal{A}_0, Y, \{R_i^{(2)}\}_{i \in \omega})$} takes for input:

(1) a predicate structure $\mathcal{A}$;

(2) a substructure $\mathcal{A}_0 \subset \mathcal{A}$, where its
universe equals to an infinite set for solutions of a formula
$\psi(x)$ in $\mathcal{A}$, the substructure generates unique
non-principal 1-type $p_{\infty}(x)\in S(\varnothing)$ and
$p_{\infty}(x)$ is realized in $\mathcal{A}_0$;

(3) an infinite set $Y$ with $Y \cap A = \varnothing$;

(4) a sequence $(R_i^{(2)})_{i\in\omega}$ of binary predicate
symbols.

We assume that $A_0$ is the domain of predicates $R_i$, $Y$ is the
range of, $\vdash R_i(x,y) \rightarrow R_0(x,y)$, $i
> 0$. The work of the operator is defined by the following schemes of formulas:

(1) $\forall x \, \exists^{\infty}y({\rm Col}_0(x) \rightarrow
R_0(x,y))$;

(2) $\forall x, x' \, (\neg (x \approx x') \rightarrow \neg
\exists y (R_0(x,y) \wedge R_0(x',y)))$, i.~e., $R_0$-images of
distinct element satisfying $\psi(x)$ are disjoint and an
equivalence relation on $Y$ with infinitely many infinite classes
is refined by the formula $R_0(x,y)$;

(3) $\forall x ({\rm Col}_n(x) \rightarrow \exists^{\infty}y
(R_0(x,y) \wedge \bigwedge\limits_{i = 1}^{n} R^{\delta_i}_i(x,y))
\wedge \neg \exists z \bigvee\limits_{i > n} R_i(x,z))$ for all
possible binary tuples $(\delta_1,\ldots, \delta_n)$, i.~e., for
any element $a \in A_0$ of color $n$, the set of solutions for the
formula $R_0(a,y)$ is divided, by $R_n(x,y)$, into $2^n$ disjoint
sets, each of which is infinite.

Thus, the set of solutions for the formula $R_0(a,y)$, where $a
\models p_{\infty}(x)$, is divided, by $R_n(x,y)$, into continuum
many disjoint sets similar Example 5. For output of the operator,
we obtain a structure $\mathcal{B}$ with continuum many
non-principal types $\{R^{\delta_i}_i(a,y) \mid i \in \omega
\setminus \{0\}\}$, and there are no prime models over the type
$p_{\infty}(x)$.

\medskip

The \emph{operator of allocation for a countable
subset}\index{Operator!of allocation for a countable subset} ${\rm
css} (\mathcal{A}, \mathbf{q}_{\omega}, \mathcal{A}_0,
\{R^{(2)}_{j}\}_{j \in \omega})$\index{${\rm css} (\mathcal{A},
\mathbf{q}_{\omega}, \mathcal{A}_0, \{R^{(2)}_{j}\}_{j \in
\omega})$} takes for input:

(1) a predicate structure $\mathcal{A}$ with a continual set
$\mathbf{q}$ of non-principal $1$-types;

(2) a countable subset $\mathbf{q}_{\omega} \subset \mathbf{q}$;

(3) a substructure $\mathcal{A}_0 \subset \mathcal{A}$ with unique
non-principal 1-type $p_{\infty}(x)\in S(\varnothing)$ and such
that $p_{\infty}(x)$ is realized in $\mathcal{A}_0$;

(4) a sequence $(R_j^{(2)})_{j\in\omega}$ of binary predicate
symbols.

Denote by ${\rm Col}_{ij}(x)$ approximations of types $q_j(x) \in
\mathbf{q}_{\omega}$, $j \in \omega$. Then the type $q_j$ is
isolated by the set of formulas $\{\neg {\rm Col}_{ij}(x) \mid i
\in \omega\}$. At the operator's work, we assume that $A_0$ is the
domain of predicates $R_{ij}$ and their range contains the set of
realizations for types in $\mathbf{q}_{\omega}$. The work of the
operator is defined by the following schemes of formulas:

(1) $\forall x ({\rm Col}_i(x) \rightarrow \bigwedge\limits_{k
\geq i} \exists^{\infty} y (R_j(x,y) \wedge {\rm Col}_{kj}(y))
\wedge \bigwedge\limits_{k < i} \neg \exists y (R_j(x,y) \wedge
{\rm Col}_{kj}(y)))$, i.~e., for any element $a \in A_0$ of $i$-th
color, there are infinitely many images of each color $k$, $k \geq
i$, and there are no images of colors $k$, $k < i$;

(2) $\forall x, x'(\neg (x \approx x') \rightarrow \neg \exists y
(R_{j}(x,y) \wedge R_{j}(x', y)))$, i.~e., images of distinct
elements belonging to $A_0$ are disjoint.

If the continual set $\mathbf{q}$ of non-principal types is
obtained by the operator ${\rm icp}$ (and there are no prime
models over each type in $\mathbf{q}$) then after passing all
colors ${\rm Col}$ by all predicates $R_j$, the countable subset
$\mathbf{q}_{\omega}$ is selected and, using a generic
construction for a structure with required properties, there
exists a prime model $\mathcal{M}_{p_{\infty}}$ over a realization
of $p_{\infty}$ and realizing exactly all types in
$\mathbf{q}_{\omega}$. If $\mathbf{q}_{\omega}$ is dense in
$\mathbf{q}$ with respect to natural topology then, assuming that
types in $\mathbf{q}_{\omega}$ are free (are not linked with $a
\in A_0$), we can remove elements in $\mathbf{q}_{\omega}$ and
obtain new prime model
$\mathcal{M}_{p_{\infty},\tilde{\mathbf{q}}_{\omega}}$,
$\tilde{\mathbf{q}}_{\omega} \subset \mathbf{q}_{\omega}$, being
an elementary submodel of
$\mathcal{M}_{p_{\infty},\mathbf{q}_{\omega}}$. But having links
of the dense set $\mathbf{q}_{\omega}$ with the type $p_{\infty}$
by predicates, the removing of a type in $\mathbf{q}_{\omega}$
leads to the removing of $p_{\infty}$. Whence applying the
operator ${\rm css}$ with input parameters, satisfying the
conditions above, there are no other (non-isomorphic) prime models
being an elementary submodel of
$\mathcal{M}_{p_{\infty},\mathbf{q}_{\omega}}$. Thus if we focus
on this property, the given operator is called the \emph{operator
of ban for downward movement}\index{Operator!of ban for downward
movement} and it is denoted by ${\rm bd}$\index{${\rm bd}
(\mathcal{A}, \mathbf{q}_{\omega}, \mathcal{A}_0,
\{R^{(2)}_{j}\}_{j \in \omega})$} with the same input parameters.

\medskip

The \emph{operator of ban for upward movement}\index{Operator!of
ban for upward movement} ${\rm bu} (\mathcal{A}, \mathcal{A}_1,
\mathcal{A}_2, Z, \{R^{(3)}_n\}_{n\in\omega})$\index{${\rm bu}
(\mathcal{A}, \mathcal{A}_1, \mathcal{A}_2, Z,
\{R^{(3)}_n\})_{n\in\omega}$} takes for input:

(1) a predicate structure $\mathcal{A}$;

(2) two disjoint substructures $\mathcal{A}_1$ and $\mathcal{A}_2$
of $\mathcal{A}$ with unique non-principal 1-types $p_1$ and
$p_2$, being realized in $\mathcal{A}_1$ and $\mathcal{A}_2$
respectively;

(3) an infinite set $Z$ such that $A_1 \cap Z = \varnothing$ and
$A_2 \cap Z = \varnothing$;

(4) a sequence $(R_n^{(3)})_{n\in\omega}$ of ternary predicate
symbols.

We denote approximations of $p_1$ and $p_2$ by ${\rm Col}_{i1}$
and ${\rm Col}_{i2}$, $i \in \omega$, respectively. The set $A_1
\times A_2$ is the domain of predicates $R_{n}$, and $Z$ is their
range, $\vdash R_n(x,y,z) \rightarrow R_0(x,y,z)$, $i
> 0$. The work of the
operator is defined by the following schemes of formulas:

(1) $\forall x, y({\rm Col}_{01}(x) \wedge {\rm Col}_{02}(y)
\rightarrow \exists^{\infty} z R_0(x, y, z))$;

(2) $\forall x, y, x', y'(\neg( x \approx x') \wedge \neg (y
\approx y') \rightarrow \neg \exists z (R_{n}(x,y,z) \wedge
R_{n}(x',y',z)))$, i.~e., $R_{n}$-images of distinct pairs
$(a_1,a_2) \in A_1 \times A_2$ are disjoint and the set $Z$ is
divided into infinitely many infinite equivalence classes;

(3) $\forall x, y ({\rm Col}_{k1}(x) \wedge {\rm Col}_{n2}(y)
\rightarrow \exists^{\infty}z (R_0(x,y,z) \wedge
\bigwedge\limits_{1 \leq i \leq \min(k,n)} R^{\delta_i}_i(x,y,z))
\wedge \neg \exists z \bigvee\limits_{i > \min(k,n)} R_i(x,y,z))$
for all possible binary tuples $(\delta_1,\ldots,
\delta_{\min(k,n)})$.

Hence, if a pair $(a_1,a_2)$ has the $(\infty,\infty)$-color, the
set of solutions for the formula $R_0(a_1,a_2,z)$ is divided on
continuum many parts. Thus, there is a prime model over each
realization of $p_1(x)$ and of $p_2(y)$, but there are no prime
models over types $q(x,y) \supset p_1(x) \cup p_2(y)$.

\medskip

The \emph{operator for construction of limit models over a
type}\index{Operator!for construction of limit models!over a
type}, ${\rm lmt}(p, \lambda,
\{R^{(2)}_i\}_{i\in\omega})$\index{${\rm lmt}(p, \lambda,
\{R^{(2)}_i\}_{i\in\omega})$} takes for input:

(1) a non-principal 1-type $p(x)$;

(2) a number $\lambda \in \omega+1$ of limit models over $p(x)$;

(3) a sequence $(R_i^{(2)})_{i\in\omega}$ of binary predicate
symbols.

We assume that predicates $R_i$ act on a set of realizations of
$p(x)$ such that $R_i(a, y) \vdash p(y)$ and $\models\exists y
R_i(a, y)$ and realizations $R_i(a, y)$ do not semi-isolate $a$,
where $a \models p(x)$. We construct a tree of $R_i$-extensions
over a realization $a$ of $p$. Consider sequences $i_0, \ldots,
i_n, \ldots \in 2^\omega$ correspondent to pathes $R_{i_0}(a, a_1)
\wedge \ldots \wedge R_{i_n}(a_n, a_{n+1}) \wedge \ldots$. There
are $2^{\omega}$ extensions. As shown in \cite{SuLP, Su072}, given
number $\lambda \in \omega + 1$ of limit models can be obtained by
some family of identities.

For $n\in\omega\setminus\{0\}$ limit models, we use the following
identities:

(1) $n - 1 \approx m$, $m \geq n$,

(2) $mm \approx m$, $m < n$,

(3) $n_1 n_2 \ldots n_s \approx n_s$, $\min \{n_1, n_2, \ldots,
n_{s-1}\} > n_s$.

For countably many limit models, we introduce identities:

(1) $nn \approx n$, $n \in \omega$,

(2) $n_1 n_2 \ldots n_s \approx n_s$, $\min \{n_1, n_2, \ldots,
n_{s-1}\} > n_s$,

(3) $n_1 n_2 \approx n_1 (n_1+1) (n_2 +2) \ldots (n_2 - 1)n_2$,
$n_1 < n_2$.

\medskip

The \emph{operator for construction of limit models over a
$\leq_{RK}$-sequence}\index{Operator!for construction of limit
models!over a $\leq_{RK}$-sequence}\index{${\rm lms}((q_n)_{n \in
\omega}, \lambda, \{R^{(2)}_i\}_{i \in \omega})$} $${\rm
lms}((q_n)_{n \in \omega}, \lambda, \{R^{(2)}_i\}_{i \in
\omega})$$ takes for input:

(1) a $\leq_{RK}$-sequence $(q_n)_{n \in \omega}$;

(2) a number $\lambda \in \omega+1$ of limit models over the
sequence $(q_n)_{n \in \omega}$;

(3) a sequence $(R_i^{(2)})_{i\in\omega}$ of binary predicate
symbols.

Consider types $q_n$ and $q_{n+1}$. Since they belong to the
$\leq_{RK}$-sequence, there is a formula $\varphi(x,y)$ such that
$q_{n+1}(y) \cup \{\varphi(x,y)\}$ is consistent and $q_{n+1}(y)
\cup \{\varphi(x,y)\} \vdash q_{n}(x)$. We assume that predicates
$R_i$ act so that $R_i(x,y) \vdash \varphi(x,y)$ and for every $a
\models q_{n+1}(y)$, $R_i(x,a) \vdash q_{n}(x)$. Below we consider
numbers $i$ instead of predicates $R_i$. Then for the
$\leq_{RK}$-sequence, there are $\omega^{\omega}$ sequences $i_l,
\ldots, i_k, \ldots$ correspondent to $R_l(a_n, a_{n-1}) \wedge
\ldots \wedge R_k(a_j, a_{j-1}) \wedge \ldots $, where $a_n
\models q_{n}(x),$ $\ldots,$ $a_1 \models q_1(x)$.

By the sequence $(q_n)_{n\in\omega}$, we construct sequences of
prime models $\mathcal{M}_{q_n}$ over realizations of $q_n$, where
$(n+1)$-th model is an elementary extension of $n$-th one. Any
limit model is a union of countable chain of a sequence of prime
models over tuples. Predicates $R_i$, $i \in \omega$, link types
in $(q_n)$ leading to required number of limit models. As shown in
\cite{SuLP, Su08}, the problem of extension of a theory producing
a given number of limit models over $(q_n)$ is reduced to a
factorization of the set $\omega^{\omega}$ by an identification of
some words such that the result of this factorization contains as
many classes as there are limit models.

For $n\in\omega\setminus\{0\}$ limit models, we use the following
identities:

(1) $n - 1 \approx m$, $m \geq n$;

(2) $n_0 n_1 \ldots n_s \approx \underbrace{n_s \ldots n_s}_{s+1
\, {\textrm{times}}}$, $\max\{n_0, n_1, \ldots, n_{s-1}\} < n_s$.

For countably many limit models, we take identities:

(1) $n_0 n_1 \ldots n_s \approx \underbrace{n_s \ldots n_s}_{s+1
\, {\textrm{times}}}$, $\max\{n_0, n_1, \ldots, n_{s-1}\} < n_s$;

(2) $n_0 n_1 \ldots n_s \approx n_0(n_0+1)\ldots(n_0+s)$, $n_0 + s
\leq n_s$;

(3) $n_0 n_1 \ldots n_s \approx
n_0(n_0+1)\ldots(n_0+t)\underbrace{(n_0+t)\ldots(n_0+t)}_{s - t \,
\textrm{times}}$, $n_0 + s$, $n_0 + t = n_s$, $t > 0$, $s > t$.

\bigskip
\centerline{\bf 7. Distributions of prime and limit models}
\centerline{\bf for finite Rudin--Keisler preorders}
\medskip

If $\widetilde{{\bf M}}$ is a $\sim_{RK}$-class containing an
isomorphism type ${\bf M}$ of a prime model over a tuple, then as
usual we denote by ${\rm IL}(\widetilde{{\bf M}})$ the number of
limit models, being unions of elementary chains of models, whose
isomorphism types belong to the class $\widetilde{{\bf M}}$.

Clearly, for theories $T$ with finite structures ${\rm RK}(T)$,
any limit model is limit over a type.

The following two theorems show that for $p$-Ehrenfeucht small
theories, the number of countable models is defined by the number
of prime models over tuples and by the distribution function ${\rm
IL}$ of numbers of limit models over types. Assuming Continuum
Hypothesis, all possible basic characteristics are realized.

\medskip
{\bf THEOREM 7.1} \cite{SuLP, Su072}. {\em Any  small  theory $T$
with  a  finite \ Rudin--Keisler preorder satisfies the following
conditions:

{\rm (a)} ${\rm RK}(T)$ contains the least element ${\bf M}_0$
{\rm (}the isomorphism type of a prime model{\rm )}, and ${\rm
IL}(\widetilde{{\bf M}_0})=0$;

{\rm (b)} ${\rm RK}(T)$ contains the greatest $\sim_{\rm
RK}$-class $\widetilde{{\bf M}_1}$ {\rm (}the class of~isomorphism
types of all prime models over realizations of powerful types{\rm
)}, and $|{\rm RK}(T)|>1$ implies ${\rm IL}(\widetilde{{\bf
M}_1})\geq 1$;

{\rm (c)} if $|\widetilde{\bf M}|>1$, then ${\rm
IL}(\widetilde{\bf M})\geq 1$.

Moreover, we have the following {\sl decomposition
formula}\index{Formula!decomposition}:
$$
I(T,\omega)=|{\rm RK}(T)|+\sum_{i=0}^{|{\rm RK}(T)/\sim_{\rm
RK}|-1} {\rm IL}(\widetilde{{\bf M}_i}),
$$
where \ $\widetilde{{\bf M}_0},\ldots, \widetilde{{\bf M}_{|{\rm
RK}(T)/\sim_{\rm RK}|-1}}$ \ are all elements of the partially
ordered set ${\rm RK}(T)/\!\!\sim_{\rm RK}$ and ${\rm
IL}(\widetilde{{\bf
M}_i})\in\omega\cup\{\omega,\omega_1,2^\omega\}$ for each $i$.}

\medskip
{\bf THEOREM 7.2} \cite{SuLP, Su072}. {\em For any finite
preordered set $\langle X;\leq\rangle$ with the least element
$x_0$ and the greatest class $\widetilde{x_1}$ in the ordered
factor set $\langle X;\leq\rangle/\!\!\sim$ with respect to $\sim$
{\rm (}where $x\sim y\Leftrightarrow x\leq y\mbox{ and }y\leq
x${\rm )}, and for any function $f\mbox{\rm :
}X/\!\!\sim\:\to\omega\cup\{\omega,2^\omega\}$, satisfying the
conditions $f(\widetilde{x_0})=0$, $f(\widetilde{x_1})>0$ for
$|X|>1$, and $f(\widetilde{y})>0$ for $|\widetilde{y}|>1$, there
exist a~small theory $T$ and an isomorphism $g\mbox{\rm : }\langle
X;\leq\rangle\:\widetilde{\to}\:{\rm RK}(T)$ such that ${\rm
IL}(g(\widetilde{y}))=f(\widetilde{y})$ for any $\widetilde{y}\in
X/\!\!\sim$.}

\medskip
Note that, by criterion of existence of prime model, an unsmall
theory $T$ is $p$-categorical if and only if there is a unique
$\equiv_{\rm RK}$-class $S\subset S(T)$ such that for any
realization $\bar{a}$ of some (any) type in $S$ every consistent
formula $\varphi(\bar{x},\bar{a})$ is an ${\rm i}$-formula.

Similarly, an unsmall theory $T$ is $p$-Ehrenfeucht if and only if
there are finitely many pairwise non-$\equiv_{\rm RK}$-equivalent
types $p_j$, $j<n$, $1<n<\omega$, such that for any $j$ and for
some (any) realization $\bar{a}_j$ of $p_j$ every consistent
formula $\varphi(\bar{x},\bar{a}_j)$ is an ${\rm i}$-formula.

The proofs of the following assertions repeat according proofs for
the class of small theories \cite{SuLP, Su041, BSV}.

\medskip
{\bf PROPOSITION 7.3.} {\em If ${\cal M}_p$ and ${\cal M}_q$ are
domination-equivalent non-iso\-mor\-phic models then there exist
models that are limit over the type $p$ and over the type $q$.}

\medskip
{\bf PROPOSITION 7.4.} {\em If types $p_1$ and $p_2$ are
domination-equi\-valent,  and  there  exists  a  limit  model over
$p_1$,  then  there  exists  a~model that is limit over $p_1$ and
over $p_2$.}

\medskip
{\bf THEOREM 7.5.} {\em Let $p(\bar{x})$ be a complete type of a
countable theory $T$. The following conditions are equivalent:

$(1)$ there exists a limit model over $p$;

$(2)$ there exists a model ${\cal M}_p$ and the relation $I_p$ of
isolation on a set of realizations of $p$ in a {\rm (}any{\rm )}
model ${\cal M}\models T$ realizing $p$ is non-symmetric;

$(3)$ there exists a model ${\cal M}_p$ and, in some {\rm
(}any{\rm )} model ${\cal M}\models T$ realizing $p$, there exist
realizations $\bar{a}$ and $\bar{b}$ of $p$ such that the type
${\rm tp}(\bar{b}/\bar{a})$ is principal and $\bar{b}$ does not
semi-isolate $\bar{a}$ and, in particular, ${\rm SI}_p$ is
non-symmetric on on the set of realizations of $p$ in ${\cal M}$.}

\medskip
By Proposition 7.3, we have the following analogue of Theorem 7.1
for the class ${\cal T}_c$.

\medskip
{\bf PROPOSITION 7.6.} {\em Every theory $T\in{\cal T}_c$ with a
finite Rudin--Keisler preorder satisfies the following: if
$|\widetilde{\bf M}|>1$ then ${\rm IL}(\widetilde{\bf M})\geq 1$.
Moreover, we have the following {\sl decomposition
formula}\index{Formula!decomposition}:
$$
I(T,\omega)=|{\rm RK}(T)|+\sum_{i=0}^{|{\rm RK}(T)/\sim_{\rm
RK}|-1} {\rm IL}(\widetilde{{\bf M}_i})+{\rm NPL}(T),
$$
where \ $\widetilde{{\bf M}_0},\ldots, \widetilde{{\bf M}_{|{\rm
RK}(T)/\sim_{\rm RK}|-1}}$ \ are all elements of the partially
ordered set ${\rm RK}(T)/\!\!\sim_{\rm RK}$ and ${\rm
IL}(\widetilde{{\bf
M}_i})\in\omega\cup\{\omega,\omega_1,2^\omega\}$ for each $i$,
$0\leq{\rm NPL}(T)\leq 2^\omega$.}

\medskip
The following theorem is an analogue of Theorem 7.2 for the class
${\cal T}_c$.

\medskip
{\bf THEOREM 7.7.} {\em For any finite preordered set $\langle
X;\leq\rangle$ and for any function $f\mbox{\rm :
}X/\!\!\sim\:\to\omega\cup\{\omega,2^\omega\}$ such that
$f(\widetilde{x})>0$ for $|\widetilde{x}|>1$ {\rm (}where $x\sim
y\Leftrightarrow x\leq y\mbox{ and }y\leq x${\rm )}, there exists
a theory $T\in{\cal T}_c$ {\rm (}without prime models{\rm )} and
an isomorphism $g\mbox{\rm : }\langle
X;\leq\rangle\:\widetilde{\to}\:{\rm RK}(T)$ such that ${\rm
IL}(g(\widetilde{x}))=f(\widetilde{x})$ for any $\widetilde{x}\in
X/\!\!\sim$.}

\medskip
PROOF. Denote the cardinality of $X$ by $m$ and consider the
theory $T_0$ of unary predicates $P_i$, $i < m$, forming a
partition of a set $A$ on $m$ disjoint infinite sets with a
coloring ${\rm Col}$: $A \rightarrow \omega\cup\{\infty\}$ such
that for any $i < m$, $j \in \omega$, there are infinitely many
realizations for each type $\{ {\rm Col}_j(x) \wedge P_i(x)\}$,
$\{\neg {\rm Col}_j(x) \mid j \in \omega \} \cup \{P_i(x)\} =
p_i(x)$. In this case, each set of formulas isolates a complete
type.

Let $X_1, \ldots,  X_n$ be connected components of the preordered
set $\langle X;\leq\rangle$, consisting of $m_1, \ldots, m_n$
elements respectively, $m_1+\ldots +m_n=m$. Now we assume that
each element in $X$ corresponds to a predicate $P_i$, $i < m$.

We expand the theory $T_0$ to a theory $T_1$ by binary predicates
$Q_{kl}$, whose domain coincides with the set of solutions for the
formula $P_k(x)$ and the range is the set of solutions for the
formula $P_l(x)$; we link types $p_k$ and $p_l$ if correspondent
elements $x_k$ and $x_l$ in $X$ belong to a common connected
component and $x_l$ covers $x_k$. Moreover, the coloring ${\rm
Col}$ will be $1$-inessential and $Q_{kl}$-ordered \cite{SuLP}:

(1) for any $i \geq j$, there are elements $x,y \in M$ such that
$$\models {\rm Col}_i(x) \wedge {\rm Col}_j(y) \wedge Q_{kl}(x,y) \wedge P_{k}(x) \wedge P_{l}(y);$$

(2) if $i < j$ then there are no elements $u,v \in M$ such that
$$\models {\rm Col}_i(u) \wedge {\rm Col}_j(v) \wedge Q_{kl}(u,v) \wedge P_{k}(u) \wedge P_{l}(v).$$

Applying a generic construction we get that if $a \models
p_{l}(y)$ then the formula $Q_{kl}(x,a)$ is isolating and
$p_{l}(y) \cup Q_{kl}(x,y) \vdash p_{k}(x)$, moreover,
realizations of $p_k$ do not semi-isolate realizations of $p_l$.
Thus the set of non-principal 1-types $p_i(x)$ has a preorder
correspondent to the preorder $\leq$.

We construct, by induction, an expansion of theory $T_1$ to a
required theory~$T$.

On initial step, we expand the theory $T_1$ by binary predicates
$\{R_i^{(2)}\}_{i \in \omega}$ and apply the operator of continual
partition ${\rm icp}(\mathcal{A}, \mathcal{A} \upharpoonright P_0,
Y, \{R_i^{(2)}\}_{i \in \omega}) = \mathcal{B}$, where
$\mathcal{A}$ is a model of $T_1$. We consider an arbitrary
connected component $X_i$ and enumerate its elements so that if
$x_k
> x_l$ then $k > l$. On further $m_i$ steps, we apply the operator
of allocation for a countable subset ${\rm css} (\mathcal{B},
\mathbf{q}_{\omega}, \mathcal{A} \upharpoonright P_{l_i},
\{R^{(2)}_{j}\}_{j \in \omega})$, where $l_1, \ldots, l_i$ are
numbers of elements forming the connected component $X_i$,
$\mathbf{q}_{\omega}$ is a countable dense subset of set
$\mathbf{q}$ of 1-types for the structure $\mathcal{B}$. We
organize a similar process for all connected components in $X$.
Now for all types corresponding to elements in distinct connected
components and to maximal elements in a common component, we apply
the operator of ban for upward movement ${\rm bu} (\mathcal{A},
\mathcal{A} \upharpoonright P_{i}, \mathcal{A} \upharpoonright
P_{j},\{R^{(3)}_{\Delta}\})$, expanding the theory by disjoint
families ternary predicates $R^{(3)}_n$, $n\in\omega$.

The required number of limit models can be done by application,
for each $g(\widetilde{x})$, of the operator ${\rm
lmt}(g(\widetilde{x}), f(\widetilde{x}),
\{R_i^{g(\widetilde{x})}\}_{i \in \omega})$ expanding the theory
by predicates $R_i^{g(\widetilde{x})}$ for each
$g(\widetilde{x})$.~$\Box$

\medskip
By the proof of Theorem 7.7, positive values $P(T)$ for the class
${\cal T}_c$ can be defined by prime models, being not prime over
$\varnothing$. Modifying the proof, one can realize an arbitrary
finite preordered set $\langle X;\leq\rangle$ with the least
element by ${\rm RK}(T)$ for a theory $T\in{\cal T}_c$ with a
prime model over $\varnothing$.

By the construction for the proof of Theorem 7.7, we get

\medskip
{\bf COROLLARY 7.8.} {\em For any cardinalities
$\lambda_1\in\omega\setminus\{0\}$ and
$\lambda_2\in\omega\cup\{\omega,2^\omega\}$ there is a theory
$T\in{\cal T}_c$ such that ${\rm
cm}_3(T)=(\lambda_1,\lambda_2,2^\omega)$.}

\bigskip
\centerline{\bf 8. Distributions of prime and limit models}
\centerline{\bf for countable Rudin--Keisler preorders}
\medskip

We say \cite{SuLP, Su08} that a family ${\bf Q}$ of $\leq_{\rm
RK}$-sequences  ${\bf q}$ of types {\em represents}\index{Family
of $\leq_{\rm RK}$-sequences representing $\leq_{\rm
RK}$-sequence} a $\leq_{\rm RK}$-sequence ${\bf q}'$ of types if
any limit model over ${\bf q}'$ is limit over some ${\bf q}\in
{\bf Q}$.

\medskip
{\bf THEOREM 8.1} \cite{SuLP, Su08}. {\em Any small theory $T$
satisfies the following conditions:

{\rm (a)} the structure ${\rm RK}(T)$ is upward directed and has
the least element ${\bf M}_0$ {\rm (}the isomorphism type of prime
model of $T${\rm )}, \ ${\rm IL}(\widetilde{{\bf M}_0})=$~$0$;

{\rm (b)} if ${\bf q}$ is a $\leq_{\rm RK}$-sequence of
non-principal types $q_n$, $n\in\omega$, such that each type $q$
of $T$ is related by $q\leq_{\rm RK} q_n$ for some~$n$, then there
exists a limit model over ${\bf q}$; in particular, $I_l(T)\geq 1$
and the countable saturated model is limit over ${\bf q}$, if
${\bf q}$ exists;

{\rm (c)} if ${\bf q}$ is a $\leq_{\rm RK}$-sequence of types
$q_n$, $n\in\omega$, and $({\cal M}_{q_n})_{n\in\omega}$ is an
elementary chain such that any co-finite subchain does not consist
of pairwise isomorphic models, then there exists a limit model
over~${\bf q}$;

{\rm (d)} if ${\bf q}'=(q'_n)_{n\in\omega}$  is a subsequence of
$\leq_{\rm RK}$-sequence ${\bf q}$, then any limit model over
${\bf q}$ is limit over ${\bf q}'$;

{\rm (e)} if ${\bf q}=(q_n)_{n\in\omega}$ and ${\bf
q}'=(q'_n)_{n\in\omega}$ are $\leq_{\rm RK}$-sequences of types
such that for some $k,m\in\omega$, since some $n$, any types
$q_{k+n}$ and $q'_{m+n}$ are related by ${\cal
M}_{q_{k+n}}\simeq{\cal M}_{q'_{m+n}}$, then any model ${\cal M}$
is limit over ${\bf q}$ if and only if ${\cal M}$ is limit
over~${\bf q}'$.

Moreover, the following  {\sl decomposition
formula}\index{Formula!decomposition} holds:
$$
I(T,\omega)=|{\rm RK}(T)|+\sum\limits_{{\bf q}\in {\bf Q}}{\rm
IL}_{\bf q},
$$
where ${\rm IL}_{\bf q}\in\omega\cup\{\omega,\omega_1,2^\omega\}$
is the number of limit models related to the $\leq_{\rm
RK}$-sequence ${\bf q}$ and not related to extensions and to
restrictions of ${\bf q}$ that used for the counting of all limit
models of $T$, and the family ${\bf Q}$ of $\leq_{\rm
RK}$-sequences of types represents all $\leq_{\rm RK}$-sequences,
over which limit models exist.}

\medskip
{\bf THEOREM 8.2} \cite{SuLP, Su08}. {\em Let $\langle
X,\leq\rangle$ be at most countable upward directed preordered set
with a least element $x_0$, $f\mbox{\rm :
}Y\to\omega\cup\{\omega,2^\omega\}$ be a function with at most
countable set $Y$ of $\leq_0$-sequences, i.~e., of sequences in
$X\setminus\{x_0\}$ forming $\leq$-chains, and satisfying the
following conditions:

{\rm (a)} $f(y)\geq 1$ if for any $x\in X$ there exists some $x'$
in the sequence $y$ such that $x\leq x'$;

{\rm (b)} $f(y)\geq 1$ if any co-finite subsequence of $y$ does
not contain pairwise equal elements;

{\rm (c)} $f(y)\leq f(y')$ if $y'$ is a subsequence of $y$;

{\rm (d)} $f(y)=f(y')$ if $y=(y_n)_{n\in\omega}$ and
$y'=(y'_n)_{n\in\omega}$ are sequences such that there exist some
$k,m\in\omega$ for which $y_{k+n}=y'_{m+n}$ since some $n$.

Then there exists a small theory  $T$ and an isomorphism
$$g\mbox{\rm : }\langle X,\leq\rangle\:\widetilde{\to}\:{\rm RK}(T)$$ such that
any value $f(y)$ is equal to the number of limit models over
$\leq_{\rm RK}$-sequence $(q_n)_{n\in\omega}$, correspondent to
the $\leq_0$-sequence $y=(y_n)_{n\in\omega}$, where $g(y_n)$ is
the isomorphism type of the model ${\cal M}_{q_n}$, $n\in\omega$.}

\medskip
Repeating the proof of Theorem 8.1, we obtain

\medskip
{\bf THEOREM 8.3.} {\em Any theory $T\in\mathcal{T}_c$ satisfies
the following conditions:

{\rm (a)} if ${\bf q}$ is a $\leq_{\rm RK}$-sequence of types
$q_n$, $n\in\omega$, and $({\cal M}_{q_n})_{n\in\omega}$ is an
elementary chain such that any co-finite subchain does not consist
of pairwise isomorphic models, then there exists a limit model
over~${\bf q}$;

{\rm (b)} if ${\bf q}'=(q'_n)_{n\in\omega}$  is a subsequence of
$\leq_{\rm RK}$-sequence ${\bf q}$, then any limit model over
${\bf q}$ is limit over ${\bf q}'$;

{\rm (c)} if ${\bf q}=(q_n)_{n\in\omega}$ and ${\bf
q}'=(q'_n)_{n\in\omega}$ are $\leq_{\rm RK}$-sequences of types
such that for some $k,m\in\omega$, since some $n$, any types
$q_{k+n}$ and $q'_{m+n}$ are related by ${\cal
M}_{q_{k+n}}\simeq{\cal M}_{q'_{m+n}}$, then any model ${\cal M}$
is limit over ${\bf q}$ if and only if ${\cal M}$ is limit
over~${\bf q}'$.

Moreover, the following  {\sl decomposition
formula}\index{Formula!decomposition} holds:
$$
I(T,\omega)=|{\rm RK}(T)|+\sum\limits_{{\bf q}\in {\bf Q}}{\rm
IL}_{\bf q}+{\rm NPL}(T),
$$
where ${\rm IL}_{\bf q}\in\omega\cup\{\omega,\omega_1,2^\omega\}$
is the number of limit models related to the $\leq_{\rm
RK}$-sequence ${\bf q}$ and not related to extensions and to
restrictions of ${\bf q}$ that used for the counting of all limit
models of $T$, and the family ${\bf Q}$ of $\leq_{\rm
RK}$-sequences of types represents all $\leq_{\rm RK}$-sequences,
over which limit models exist.}

\medskip
Similarly Theorem 7.2, Theorem 8.2 has a generalization for the
class $\mathcal{T}_c$:

\medskip
{\bf THEOREM 8.4.} {\em Let $\langle X,\leq\rangle$ be at most
countable preordered set, $f\mbox{\rm :
}Y\to\omega\cup\{\omega,2^\omega\}$ be a function with at most
countable set $Y$ of $\leq$-sequences, i.~e., of sequences in $X$
forming $\leq$-chains, and satisfying the following conditions:

{\rm (a)} $f(y)\geq 1$ if any co-finite subsequence of $y$ does
not contain pairwise equal elements;

{\rm (b)} $f(y)\leq f(y')$ if $y'$ is a subsequence of $y$;

{\rm (c)} $f(y)=f(y')$ if $y=(y_n)_{n\in\omega}$ and
$y'=(y'_n)_{n\in\omega}$ are sequences such that there exist some
$k,m\in\omega$ for which $y_{k+n}=y'_{m+n}$ since some $n$.

Then there exists a theory $T\in\mathcal{T}_c$ and an isomorphism
$$g\mbox{\rm : }\langle X,\leq\rangle\:\widetilde{\to}\:{\rm RK}(T)$$ such that
any value $f(y)$ is equal to the number of limit models over
$\leq_{\rm RK}$-sequence $(q_n)_{n\in\omega}$, correspondent to
the $\leq$-sequence $y=(y_n)_{n\in\omega}$, where $g(y_n)$ is the
isomorphism type of the model ${\cal M}_{q_n}$, $n\in\omega$.}

\medskip
PROOF. We assume that $X$ is countable since for finite $X$, the
proof repeats the construction for the proof of Theorem 7.7. Now
we consider the theory $T_0$ of unary predicates $P_i$, $i \in
\omega$, forming, with the type $p_{\infty}(x) = \{\neg P_i(x)
\mid i \in \omega \}$, a partition of a set $A$ by disjoint
infinite classes with a coloring ${\rm Col}$: $A \rightarrow
\omega\cup\{\infty\}$ such that for any $i,j \in \omega$, there
are infinitely many realizations for each of types $\{ {\rm
Col}_j(x) \wedge P_i(x)\}$, $\{\neg {\rm Col}_j(x) \mid i \in
\omega \} \cup \{P_i(x)\} = p_i(x)$, $\{ {\rm Col}_j(x) \} \cup
p_{\infty}(x)$, $\{\neg {\rm Col}_j(x) \mid j \in \omega \}\cup
p_{\infty}(x)$. Here, each set of formulas isolates a complete
type. We link the type $\{\neg {\rm Col}_j(x) \mid j \in
\omega\}\cup p_{\infty}(x)$ with the type $p_0(x)$ by an extension
of $T_0$ to a theory $T_1$ with a binary predicate $Q_0$ such that
for all $j \in \omega$, we have:

(1) $\forall x, y \left( {\rm Col}_j(x) \wedge P_0(x) \wedge
Q_0(x,y) \rightarrow {\rm Col}_j(y) \wedge P_j(y)\right);$

(2) $\forall x, y \left( {\rm Col}_j(y) \wedge P_j(y) \wedge
Q_0(x,y) \rightarrow {\rm Col}_j(x) \wedge P_0(x)\right);$

(3) $Q_0$ is a bijection between sets of solutions for the
formulas ${\rm Col}_j(x) \wedge P_0(x)$ and ${\rm Col}_j(y) \wedge
P_j(y)$.

These conditions allow not to care about the type $p_{\infty}(x)$
with respect to the existence of prime model over it, since
$p_0(x)$ and $p_{\infty}(x)$ are strongly ${\rm RK}$-equivalent.

Let $X_1, \ldots, X_n,\ldots$ be connected components in the
preordered set $\langle X,\leq\rangle$. We consider a one-to-one
correspondence between $X$ and the set of predicates $P_i(x)$, $i
\in \omega$.

Similar the proof of Theorem 7.7, we expand the theory $T_1$ to a
theory $T_2$ by binary predicates $Q_{kl}$ with domains $P_k(x)$
and ranges $P_l$, and link types $p_k$ and $p_n$ if correspondent
elements in $X$ lay in common connected component and an element
$x_l$ corresponding to $p_l$ covers an element $x_k$ corresponding
to $p_k$. Moreover, using a generic construction, the coloring
${\rm Col}$ should be $1$-inessential and $Q_{kl}$-ordered.

The further proof repeats arguments for the proof of Theorem 7.7,
where the operator ${\rm css}$ of allocation for a countable set
is applied countably many times, for non-principal types
corresponding to elements in $X$. In this case, if non-principal
types are not exhausted, we apply the operator ${\rm icp}$ of
continual partition for remaining types.

For the required number of limit models with respect to a sequence
$(q_n)_{n \in \omega}$, we expand the theory by predicates
$R^{(q_n)}_i$, $i \in \omega$, and apply the operator ${\rm
lms}((q_n)_{n \in \omega}, f(y), \{R^{(q_n)}_i\}_{i \in \omega}),$
where $y$ is a sequence in $Y$ correspondent to the sequence
$(q_n)_{n \in \omega}$.~$\Box$

\medskip
By the construction for the proof of Theorem 8.4, we obtain

\medskip
{\bf COROLLARY 8.5.} {\em For any cardinality
$\lambda\in\omega\cup\{\omega,2^\omega\}$, there is a theory
$T\in{\cal T}_c$ such that ${\rm
cm}_3(T)=(\omega,\lambda,2^\omega)$.}

\bigskip
\centerline{\bf 9. Interrelation of classes ${\bf P}$, ${\bf L}$,
and ${\bf NPL}$} \centerline{\bf in theories with continuum many
types.} \centerline{\bf Distributions of triples ${\rm cm}_3(T)$
in the class $\mathcal{T}_c$}

\medskip
{\bf THEOREM 9.1.} {\em Let $\langle X,\leq\rangle$ be at most
countable preordered set, where $X$ is a disjunctive union of some
sets $P$ and ${\rm NPL}$, $f\mbox{\rm :
}Y\to\omega\cup\{\omega,2^\omega\}$ be a function with at most
countable set $Y$ of $(P,\leq)$-sequences, i.~e., of sequences in
$P$ forming $\leq$-chains, and satisfying the following
conditions:

{\rm (a)} $f(y)\geq 1$ if any co-finite subsequence of $y$ does
not contain pairwise equal elements;

{\rm (b)} $f(y)\leq f(y')$ if $y'$ is a subsequence of $y$;

{\rm (c)} $f(y)=f(y')$ if $y=(y_n)_{n\in\omega}$ and
$y'=(y'_n)_{n\in\omega}$ are sequences such that there exist some
$k,m\in\omega$ for which $y_{k+n}=y'_{m+n}$ since some $n$.

Then there is a theory $T\in\mathcal{T}_c$ and an isomorphism
$g\mbox{\rm : }\langle X,\leq\rangle\:\widetilde{\to}\:{\rm
CM}_0(T)$ to a substructure ${\rm CM}_0(T)=\langle{\bf
CM}_0(T);\leq_{\rm RK}\rangle$ of ${\rm CM}(T)$, with ${\bf
CM}_0(T)\subset{\bf P}(T)\cup{\bf NPL}(T)$ and satisfying the
following:

{\rm (1)} $g(P)={\bf P}(T)$, $g({\rm NPL})={\bf CM}_0(T)\cap{\bf
NPL}(T)$;

{\rm (2)} each value $f(y)$ is equal to the number of limit models
over a $\leq_{\rm RK}$-sequence $(q_n)_{n\in\omega}$ correspondent
to the $\leq$-sequence $y=(y_n)_{n\in\omega}$, where $g(y_n)$ is
the isomorphism type of the model ${\cal M}_{q_n}$, $n\in\omega$.}

\medskip
PROOF. The construction of preordered set of types, isomorphic to
the structure $\langle X,\leq\rangle$ and without prime models
over the type $p_0$ is similar the proof of Theorem  8.4. Then for
each non-principal type $p_i$, correspondent to an element in $P$,
we apply the operator of allocation for a countable subset ${\rm
dss}(\mathcal{A}, \textbf{q}_{\omega}, \mathcal{A} \upharpoonright
P_i, \{R_n\}_{n \in \omega})$. If there are types $p_i$,
correspondent to elements in $NPL$, we apply, for these types, the
operator of continual partition ${\rm icp}(\mathcal{A},
\mathcal{A} \upharpoonright P_i, Z, \{R_n\}_{n \in \omega})$. For
all types, corresponding to elements in distinct connected
components in $\langle X,\leq\rangle$ as well as to maximal
elements in a common component, we apply the operator of ban for
upward movement. For the removing of prime models over remaining
continuum many types, we apply, for $n$-tuples of elements the
operator of continual partition, using $(n+1)$-ary predicates. The
required number of limit models is obtained by the operator for
construction of limit models over a sequence of types.~$\Box$

\medskip
{\bf THEOREM 9.2.} {\em In the conditions of Theorem $9.1$, there
is a theory $T\in\mathcal{T}_c$ and an isomorphism $g\mbox{\rm :
}\langle X,\leq\rangle\,\,\widetilde{\to}\,\,{\rm CM}_0(T)$ to a
substructure $${\rm CM}_0(T)=\langle{\bf CM}_0(T);\leq_{\rm
RK}\rangle$$ of ${\rm CM}(T)$, with ${\bf CM}_0(T)\subset{\bf
P}(T)\cup{\bf NPL}(T)$ and satisfying the following:

{\rm (1)} $g(P)={\bf CM}_0(T)\cap{\bf P}(T)$, $g({\rm NPL})={\bf
NPL}(T)$;

{\rm (2)} each value $f(y)$ is equal to the number of limit models
over a $\leq_{\rm RK}$-sequence $(q_n)_{n\in\omega}$ correspondent
to the $\leq$-sequence $y=(y_n)_{n\in\omega}$, where $g(y_n)$ is
the isomorphism type of the model ${\cal M}_{q_n}$, $n\in\omega$.}

\medskip
PROOF is similar the proof of Theorem 9.1 with the only difference
that before we use the operator of continual partition and then,
if non-principal types $p_i$ are not exhausted, we apply the
operator of allocation for a countable set. For getting prime
models over remaining continuum many types, we apply, for
$n$-tuples of elements the operator of allocation for a countable
set, using $(n+1)$-ary predicates. The required number of limit
models is obtained by the operator for construction of limit
models over a sequence of types.~$\Box$

\medskip
By the construction for the proof of Theorem 9.2, we obtain

\medskip
{\bf COROLLARY 9.3.} {\em For any cardinalities
$\lambda\in\omega\cup\{\omega,2^\omega\}$ there is a theory
$T\in\mathcal{T}_c$ such that ${\rm
cm}_3(T)=(2^\omega,2^\omega,\lambda)$.}

\medskip
Proposition 5.4 and Corollaries 7.8, 8.5, 9.3 imply the following
analogue of Theorem 5.1 for the class $\mathcal{T}_c$.

\medskip
{\bf THEOREM 9.4.} {\em In the continuum hypothesis, for any
theory $T$ in the class $\mathcal{T}_c$ the triple ${\rm cm}_3(T)$
has one of the following values:

$(1)$ $(2^\omega,2^\omega,\lambda)$, where
$\lambda\in\omega\cup\{\omega,2^\omega\}$;

$(2)$ $(0,0,2^\omega)$;

$(3)$ $(\lambda_1,\lambda_2,2^\omega)$, where $\lambda_1\geq 1$,
$\lambda_1,\lambda_2\in\omega\cup\{\omega,2^\omega\}$.

All these values have realizations in the class $\mathcal{T}_c$.}

\medskip
In conclusion, the authors thank Evgeniy A.~Palyutin for helpful
remarks.

\end{document}